\documentclass{amsart}
\usepackage[
    doi=false,
    isbn=false,
    sortcites,
    backend=biber,
    giveninits=true,
    maxbibnames=99,
    ]{biblatex} 
\AtBeginBibliography{\setlength{\emergencystretch}{3em}}

\usepackage{amsmath, amssymb, amsthm}
\usepackage{thmtools}
\usepackage{fullpage}
\usepackage{float}
\usepackage{capt-of}
\usepackage{placeins}
\usepackage[hidelinks]{hyperref}
\usepackage{cleveref}
\usepackage{tikz}
\usetikzlibrary{fit}

\numberwithin{equation}{section}

\newcommand{\nn}{\mathcal{N}}
\newcommand{\realization}{\mathcal{R}}
\newcommand{\interpretation}{\mathcal{I}}
\newcommand{\relu}{\sigma}
\newcommand{\identity}{\rm id}
\newcommand{\order}{\mathcal{O}}

\newcommand{\Q}{\mathbb{Q}}
 
\newcommand{\fz}[1]{#1}
\newcommand{\mf}[1]{#1}
\newcommand{\bg}[1]{#1}
\def\AA{\mathcal{ A}}

\def\R{{\mathbb R}}
\def\E{{\mathbb E}}
\def\N{{\mathbb N}}
\def\Z{{\mathbb Z}}

\def\NN{{\mathcal N}}

\def\RR{{\mathcal R}}

\def\TT{{\mathcal T}}

\newcommand{\cA}{\mathcal{A}}
\newcommand{\cR}{\mathcal{R}}

\def\be{{\boldsymbol{e}}}
\def\bv{{\boldsymbol{v}}}

\def\ba{{\boldsymbol{a}}}
\def\bb{{\boldsymbol{b}}}

\def\bW{{\boldsymbol{W}}}
\def\bb{{\boldsymbol{b}}}

\def\norm#1#2{\|#1\|_{#2}}

\def\set#1#2{\big\{#1\,:\,#2\big\}}

\def\eps{\varepsilon}
\def\normL2#1#2{\|#1\|_{L^2(#2)}}

\newtheorem{theorem}{Theorem}[section]
\newtheorem{proposition}[theorem]{Proposition}
\newtheorem{lemma}[theorem]{Lemma}
\newtheorem{corollary}[theorem]{Corollary}

\newenvironment{remark}[1][]{%
    \medskip
    \noindent
    \textbf{Remark}%
    \if\relax\detokenize{#1}\relax
    \else\ \textnormal{(#1)}%
    \fi
    \textbf{.}\itshape
}{%
    \qed
    \smallskip
}

\DeclareRobustCommand{\CSrange}[1]%
{\ensuremath{[#1]}}

\numberwithin{theorem}{section}

\title{Computational Math with Neural Networks is Hard}
\author{Michael Feischl, Benedikt Gr\"a\ss{}le and Fabian Zehetgruber}
\thanks{Funded by the Deutsche Forschungsgemeinschaft (DFG, German Research Foundation) -- Project-ID 258734477 -- SFB 1173, the Austrian Science Fund (FWF) under the special research program Taming complexity in PDE systems (grant SFB F65) as well as project I6667-N. Funding was received also from the European Research Council (ERC) under the European Union’s Horizon 2020 research and innovation programme (Grant agreement No. 101125225) and from the OeAD under the Marietta Blau-Grant.}

\begin{document}

\begin{abstract}
We establish asymptotic complexity barriers for three basic numerical tasks with neural-network data: integration of a network realization, solution of a Poisson problem with a neural-network right-hand side, and multiplication by a neural-network encoded matrix.
A common reduction encodes Boolean formulas as ReLU networks and reduces Boolean satisfiability to each task.
Any algorithm improving on the resulting convergence-rate bounds would contradict the Strong Exponential Time Hypothesis (SETH).
These bounds hold even with full access to the network representation, including all weights.
The core reductions use networks with three hidden layers, width linear in the number of Boolean variables, bounded realizations, and sparse connectivity.
A finite-precision Monte Carlo quadrature attains the optimal rate for integration up to logarithmic factors, so the lower bound for integration is sharp.
We show slightly weaker lower bounds for networks with two hidden layers and provide fast quadrature algorithms for one-hidden-layer networks.
\end{abstract}

\maketitle 

\section{Introduction}
\noindent

Neural networks describe expressive nonlinear ansatz functions and have been shown to provide excellent surrogates for high-dimensional functions in many settings.
They are competitive with specialized approximation methods such as polynomials, rational functions, sparse-grids, and tensor networks to only name few, with prominent examples in~\cite{jentzen_approx,schwab_approx,schwab_approx2,schwab_approx3}. 
Strong results are also available for more involved applications like partial differential equations~\cite{deep_ritz0,deep_ritz1,ZANG2020109409} and inverse problems~\cite{schwab_inverse,berg_inverse}.
This also includes settings that are challenging for classical approximation methods such as high-dimensional problems~\cite{schwab_approx2}, fractals~\cite{fractals}, and stochastic processes~\cite{jentzen_approx}.
Thus, two natural questions arise.
First, can we efficiently find such networks for a given problem?
Second, once they are available, can we efficiently perform computations with them?
After all, a surrogate is usually computed with the intention of using it in a subsequent algorithm.
We refer to~\cite{optimalnn} for a recent approach to the first question and study complexity barriers for the second question in this paper.
To this end, we relate these numerical tasks to the Boolean satisfiability problem ($k$-SAT).

Under the \emph{Strong Exponential Time Hypothesis} (SETH), we show for three fundamental tasks from computational mathematics that, even with full knowledge of the neural network representation, no algorithm can attain a uniformly better algebraic convergence rate than the bounds established below.
SETH is a standard hypothesis underlying many fine-grained complexity bounds and states, informally, that brute-force search is asymptotically optimal for $k$-SAT as \(k\to\infty\).
Accordingly, any algorithm violating one of these bounds would yield a $k$-SAT algorithm whose running time contradicts SETH.
Such a refutation of SETH would be unexpected and imply surprisingly fast algorithms for several other difficult problems, including NP-hard problems~\cite{ABDN:MoreConsequencesFalsifying2018,CDL+:ProblemsHardCNFSAT2016}.

For all three tasks, the proof is based on the same reduction argument.
Given a sparse $k$-SAT formula, we construct in polynomial time a ReLU network whose realization is identically zero if the formula is unsatisfiable, whereas in the satisfiable case the relevant numerical quantity is separated from zero by an explicitly controlled gap.
A sufficiently accurate numerical algorithm that resolves this gap therefore decides the satisfiability of the original $k$-SAT formula.
Hence if the numerical algorithm was sufficiently fast, the resulting \(k\)-SAT algorithm would contradict SETH.
We give a short overview of the three tasks.

\emph{Quadrature: }
We regard the evaluation of the map
\begin{align*}
\Phi \mapsto \int_{\Omega} \realization_\Phi(x) \, dx
\end{align*}
for a neural network \(\Phi\), its realization \(\realization_\Phi\), and a domain \(\Omega\subseteq\R^d\) as a fundamental computational task.
Quadrature is a basic building block of many algorithms, including neural-network training.
For example, training PINNs typically requires approximating integral loss functions~\cite{deepRitz_E2018}, while Variational Monte Carlo involves integral-type scalar products~\cite{vmc,Gerard2022GoldstandardST}.
Related integration problems arise for the Deep Ritz method~\cite{deepRitz_E2018}, eigenvalue solvers~\cite{HAN2020109792}, and weak adversarial methods~\cite{ZANG2020109409}.
These applications require some form of quadrature that is often realised with Monte Carlo methods in practice.

Monte Carlo quadrature has a root mean-square error of order one half in the number of samples, which raises the question whether faster methods exist. Additional regularity \mf{can imply} higher convergence rates for Gaussian quadrature, quasi-Monte Carlo, or sparse-grid quadrature.
Our interest is in convergence guarantees that do not impose such additional regularity but rely only on the network representation.
Under SETH, we show that no uniformly higher algebraic convergence rate is possible \mf{already for very simple} network classes considered below.
We complement the lower bounds by a finite-bit Monte Carlo algorithm whose error has the same \mf{$L^2$-scaling as standard Monte Carlo} and attains every algebraic exponent below one half in the total runtime; see \Cref{sec:random_bit_mc}.

\emph{Solving PDEs:}
A similar question arises for the approximation of PDE solutions.
The approximation of smooth solution maps by neural networks forms the basis of many operator-learning approaches.
We refer to the overview articles~\cite{kovachki2024operatorlearningalgorithmsanalysis,azizzadenesheli2024neuraloperatorsacceleratingscientific} and the references therein, and to~\cite{neuralops,MarcatiChS24,schwab2023deepoperatornetworkapproximation,Lan2022DeepONets} for approximation-rate results.

Here we show that such regularity is essential.
Already for the linear Poisson problem $-\Delta u=f$ with Dirichlet boundary conditions, a neural-network encoded right-hand side $f=\realization_\Phi$ rules out uniformly higher algebraic convergence rates for approximating the solution map
\begin{align*}
\Phi \mapsto u_\Phi\quad\text{ with } -\Delta u_\Phi = \realization_\Phi \text{ and } u_\Phi =0 \text{ on }\partial\Omega.
\end{align*}
Our complexity barriers concern the uniform computation of the solution from the right-hand side encoded by neural networks.
These barriers do not rule out accurate neural-network representations of individual solutions or of the solution operator. They show that the existence of such representations does not by itself provide an efficient procedure for computing them from the data.

\emph{Matrix-vector multiplication:}
Nonlinear representations of high-dimensional objects have gained significant interest, with low-rank tensor formats providing a prominent example~\cite{tt1,tt2}.
In these formats, high-dimensional objects such as PDE solutions~\cite{qtt} or large matrices~\cite{ttmatrix} are stored in compressed form while still permitting efficient arithmetic.
The complexity of basic operations for other succinct matrix representations, including multi-terminal decision diagrams, has been studied in~\cite{lohrey2017succinct}.
We establish an analogous complexity barrier for matrices encoded by
neural networks.
To this end, we consider large matrices as one of the simplest
objects that can be used to store high-dimensional data. Concretely, we consider matrices $M_\Phi\in \R^{2^d\times 2^d}$ defined by
\begin{align*}
(M_{\Phi})_{ij} := \realization_\Phi(b(i)_1,\ldots,b(i)_d,b(j)_1,\ldots,b(j)_d),
\end{align*}
where $b(i)=(b(i)_1,\dots,b(i)_d)\in \{0,1\}^d$ is the binary representation of $i\in \{0,\dots,2^d-1\}$.
We prove a corresponding complexity barrier for algorithms computing such matrix-vector products.

\subsection{Related Work}
The works~\cite{liu2024neuralnetworksintegrable,RIVERA2022114710}
investigate quadrature methods for neural-network integrands, while
\cite{adaptivequad} considers adaptive quadrature and provides numerical
evidence. Quasi-Monte Carlo sampling has also been studied for the training of
neural PDE solvers~\cite{chen2021qmc}.
It is known from~\cite{t2p,grohs2025theorytopracticegapneuralnetworks} that there are no higher-order quadrature algorithms based solely on point evaluation for neural network integrands.
In our setting, the algorithm may access the complete finite representation of the network, including its weights, and the lower bounds follow from reductions from Boolean satisfiability.

We refer to~\cite{nphard,np1,np2} and the references therein for NP-hardness results for neural-network training, which include shallow ReLU networks.
Neural network verification is NP-complete in general~\cite{katz2017reluplex}, and NP-hardness persists for one-hidden-layer networks with strongly restricted parameter sets~\cite{salzer2022reachability}.
These results concern network training or the existence of inputs satisfying prescribed constraints, whereas our results concern numerical operations on a given network realization.
In particular, our lower bounds are independent of network training.

\subsection{Open Question}
The negative results of this work leave one very interesting open question: 
\begin{quote}
Is there a \emph{relevant} subset of neural networks for which higher-order algorithms exist?
\end{quote}
The term \emph{relevant} requires more explanation: the neural networks in this subset should remain expressive beyond the classical smoothness regimes for which efficient approximation methods are already available.
Our results show that sparsity, boundedness, and low-precision weights alone do not characterize such a subset.
Even for networks with two hidden layers, we obtain some lower bounds, but they are weaker than for three hidden layers.
In contrast, networks with one hidden layer can be integrated \fz{efficiently}.
\subsection{Outline of the paper}
Section~\ref{sec:assumptions} introduces the complexity-theoretic framework and the construction that encodes Boolean formulas by neural networks.
In Section~\ref{sec:quad}, we derive the quadrature lower bounds\fz{ and the matching algebraic upper bound for finite-bit Monte Carlo}.
Section~\ref{sec:computational_tasks} treats the remaining two computational tasks: approximation of Poisson solutions and matrix-vector multiplication.
Section~\ref{sec:extensions} discusses extensions and limitations, establishes fast quadrature for one-hidden-layer networks, and presents quasi-Monte Carlo experiments.

\subsection{Notation}
We write $\N:=\{1,2,3,\ldots\}$, $\N_0:=\N\cup\{0\}$, $\Q_{>0}:=\Q\cap(0,\infty)$, and
$\CSrange d:=\{1,\ldots,d\}$ for $d\in\N$, with $\CSrange0:=\varnothing$.
The cardinality of a finite set $A$ is denoted by $\#A$, while $|E|$ denotes the Lebesgue measure of a measurable set $E$.
For scalars, $|\cdot|$ denotes the absolute value and, for vectors, $|\cdot|_2=\norm{\cdot}{2}$ denotes the Euclidean norm. \fz{The notation $\operatorname{poly}(n)$ denotes a polynomial factor in $n$ whose degree is independent of $n$.}
For matrices, $\norm{\cdot}{2}$ denotes the spectral norm and $\norm{\cdot}{F}$ the Frobenius norm.
The symbols $\boldsymbol 0$, $\boldsymbol{1/2}$, and $\boldsymbol 1$ denote constant vectors of the dimension determined by context.

Throughout this paper, $\sigma(t):=\max\{0,t\}$ denotes the rectified linear unit (ReLU).
A network $\Phi$ has depth $L$ if it has $L$ affine layers, and hence $L-1$ hidden layers.
Its realization is denoted by $\mathcal R_\Phi$, and its width is the maximal layer width.
For a CNF formula $\alpha$, the quantity $\#\alpha$ denotes its number of clauses and $\ell(\alpha)$ its number of literal occurrences.
A literal is represented by a pair $(i,\gamma)$ with $i\in\CSrange n$ and $\gamma\in\{\mathrm{id},\neg\}$, where $\neg(t):=1-t$ on $[0,1]$.
Finally, $N\in\N$ denotes the work parameter supplied to a numerical algorithm, whereas $n$ denotes the number of Boolean variables in a SAT instance.
\fz{We use $d$ for a general input dimension and $s$ when the input dimension is held fixed as $N$ grows; thus $d=s$ in the fixed-dimensional setting.}

\section{From Boolean formulas to neural networks}\label{sec:assumptions}

This section introduces the complexity-theoretic framework and computational model, and develops the polynomial-time reduction that encodes SAT instances as ReLU networks with controlled architecture.

\subsection{Boolean satisfiability}
A Boolean formula is a formula in propositional logic that consists of variables that can be either $1$ (TRUE) or $0$ (FALSE) and logical connectives such as $\land$ (AND), $\lor$ (OR), and $\lnot$ (NOT).
\fz{For every Boolean formula $\alpha$ with $n\in\mathbb{N}$ variables, we define its interpretation $\interpretation_{\alpha}:\{0,1\}^n\to\{0,1\}$, which maps an assignment to the truth value of $\alpha$ under that assignment.}
A formula $\alpha$ is satisfiable if there exists an assignment $x$ that resolves to true, i.e., $\mathcal{I}_\alpha(x)=1$.
An example of a Boolean formula is
\begin{align*}
    \alpha:= (x_1 \lor (x_2 \land \lnot x_3)) \lor (\lnot x_1 \land x_3).
\end{align*}
This formula is satisfied by the assignment $x_1 = 0$, $x_2 = 1$, and $x_3 = 1$.
There may be multiple assignments that satisfy the same formula.
\fz{The Boolean satisfiability problem (SAT) asks whether a given formula $\alpha$ with $n$ variables has an assignment $x\in\{0,1\}^n$ with $\interpretation_{\alpha}(x)=1$.}

%We are interested in Boolean formula with a specific structure.
A literal is a variable $x_i$ or the negation of a variable $\lnot x_i$, and a clause is a disjunction of literals (e.g., $x_1\lor \lnot x_2\lor \ldots \lor x_n$).
A formula in conjunctive normal form (CNF) is a conjunction of clauses, for example
\begin{align}\label{eqn:CNF_sample}
    \alpha:= (x_1 \lor x_2) \land (\lnot x_1 \lor x_3 \lor x_4) \land (\lnot x_2 \lor \lnot x_3).
\end{align}
Inspired by a natural representation in computer programs, we identify a literal with the tuple $(i,\gamma)$, where $i \in \N$ is the variable index and $\gamma \in \{\neg, \identity\}$.
A clause $C$ is then a set of $q =\# C \in \N$ literals $C = \{\lambda_1, \ldots, \lambda_q\}$, and a CNF formula $\alpha$ is a set of $m = \#\alpha \in \N$ clauses $\alpha \equiv \{C_1, \ldots, C_m\}$.
For $\alpha$ from~\eqref{eqn:CNF_sample}, this leads to
\begin{align*}
	\alpha\equiv\{C_1,C_2,C_3\}
	\qquad\text{with}\;
	C_1:=\{(1,\identity),(2,\identity)\},
	\;
	C_2:=\{(1,\neg),(3,\identity),(4,\identity)\},
	\;
	C_3:=\{(2,\neg),(3,\neg)\}.
\end{align*}
We denote the set of all possible CNF formulas by $\mathcal{S}$ and define
\begin{align*}
    \mathcal{S}(n) := \set{\alpha \in \mathcal{S} }{ \text{all }(i,\gamma)\in C\in \alpha\text{ satisfy }
	i\in\CSrange{n} }.
\end{align*}
Furthermore, we will require the subset $\mathcal{S}_k(n)$, with at most $k$ literals per clause, namely
\begin{align*}
\mathcal{S}_k(n) := \set{\alpha \in \mathcal{S}(n)}{ \text{all }C\in \alpha\text{ satisfy } \#C\leq k}.
\end{align*}
Observe that every clause $C$ in $\alpha \in \mathcal{S}_k(n)$ contains at most $k$ literals,
out of the possible $2n$ different literals in $n$ variables. 
Hence any $\alpha\in \mathcal{S}_k(n)$ has at most
$\#\alpha\leq\sum_{j=1}^k\binom{2n}{j}$ clauses.

\begin{remark}
    The Tseytin transformation converts any Boolean formula into an equisatisfiable formula in CNF whose length is linear in that of the original formula; see, e.g.,~\cite{Tseitin1983}.
\end{remark}

\subsection{Algorithms and computational complexity}
\label{sec:algorithms}
Our complexity bounds are formulated through SETH and use standard discrete computation models that reflect the finite-bit arithmetic of modern computers.
\fz{We consider uniform algorithms acting on finite binary encodings.
The binary input length is the number of bits needed to store the complete input.
We use a word-RAM model: memory is organised into words, each storing a block of $\Theta(\log_2(2+B_{\rm in}))$ bits when the input has length $B_{\rm in}$.
Reading or writing a word, comparing words, and elementary arithmetic or bit operations on a bounded number of words have unit cost.
An integer that exceeds one word is stored in several words, and arithmetic on it is charged for all the word operations used.}
Uniformity means that the same finite program is used for all input lengths.
All scalar input data, constants, and outputs are rational numbers represented with finitely many bits.
In particular, this excludes arbitrary exact real constants and unit-cost real arithmetic.

For an input $\theta$ from a set of admissible inputs $\Theta$, the runtime of an algorithm $\AA$ counts the complete computation, including reading the input, arithmetic operations and comparisons, and writing the output.
In our numerical results, the input pair $(\Phi,N)$ consists of a rational neural network $\Phi$ and a work parameter $N\in\N$.
A runtime bound of $N$ therefore means that $\AA(\Phi,N)$ uses at most $N$ computational steps in total.

Our main results concern deterministic algorithms.
We also recall the corresponding randomized terminology, since the same deterministic reductions apply to randomized algorithms under randomized SETH.
\fz{A randomized algorithm may generate independent unbiased bits, each equal to $0$ or $1$ with probability $1/2$, at unit cost per bit.
A uniform random-bit algorithm is one such finite program, valid for every input length, whose randomness comes only from these bits.}
\fz{All randomized runtime bounds in this paper hold for every outcome of the random bits.
A prescribed time limit can be enforced by stopping the algorithm and returning a default result when that limit is reached; the stated success probability must then hold for this truncated algorithm.}

\noindent\textbf{Correctness of decision and approximation algorithms.}
A deterministic decision algorithm is correct if it returns the correct element of \(\{0,1\}\) for every valid input $\theta\in\Theta$.
A randomized decision algorithm has bounded error if
\begin{align*}
    \mathbb P\big(\AA(\theta)\text{ is correct}\big)\geq\frac23
    \qquad\text{for every }\theta\in\Theta.
\end{align*}
We say that an approximation algorithm for a given target \(\TT\colon\Theta\to\R\) has error at most \(\varepsilon\colon\Theta\to[0,\infty)\) if
\begin{align*}
    |\AA(\theta)-\TT(\theta)|\leq\varepsilon(\theta)
    \qquad\text{for every }\theta\in\Theta
\end{align*}
in the deterministic case, and if, for all valid inputs, this inequality holds with
probability at least \(2/3\) in the randomized
case. The reductions employ the following elementary consequence.
Suppose \(\Theta=\Theta_0\cup\Theta_1\) satisfies
\[
    \TT(\theta)\leq a\quad(\theta\in\Theta_0),
    \qquad
    \TT(\theta)\geq b\quad(\theta\in\Theta_1)
\]
for known rational numbers \(a<b\).
Then an approximation algorithm for $\TT$ with error strictly smaller than $(b-a)/2$ defines a decision algorithm that distinguishes $\Theta_0$ from $\Theta_1$ by comparing its output with $(a+b)/2$.
For randomized algorithms, the resulting decision is correct with at least the same success probability.
Repetition and the clustering argument in Lemma~\ref{lem:cluster} can amplify this probability.

\begin{lemma}[Probability amplification]\label{lem:cluster}
Let \(Y\) be a real-valued random variable and suppose that
\[
    \mathbb P(|Y-X|\leq\varepsilon)\geq\frac23
\]
for some \(X\in\R\) and \(\varepsilon>0\). Given \(\ell\) independent
samples \(y_1,\ldots,y_\ell\) of \(Y\), one can compute in
\fz{\(\mathcal O(\ell\log(\ell+1))\) comparisons and arithmetic operations} a value \(\widehat X\) such
that
\[
    \mathbb P\bigl(|\widehat X-X|\leq3\varepsilon\bigr)
    \geq
    1-\exp(-\ell/18).
\]
\end{lemma}
\begin{proof}
Let \(I:=\{i:|y_i-X|\leq\varepsilon\}\). Hoeffding's inequality gives
\[
    \mathbb P(|I|\leq\ell/2)
    \leq
    \exp\left(-2\ell(2/3-1/2)^2\right)
    =
    \exp(-\ell/18).
\]
Sort the samples and use a sliding window to find a largest subset
\(I_\star\) whose diameter is at most \(2\varepsilon\). On the event
\(|I|>\ell/2\), the set \(I\) is such a subset. Hence
\(|I_\star|>\ell/2\), so \(I_\star\cap I\neq\varnothing\). For any
\(i\in I_\star\) and \(i_0\in I_\star\cap I\),
\[
    |y_i-X|
    \leq
    |y_i-y_{i_0}|+|y_{i_0}-X|
    \leq
    3\varepsilon.
\]
We may therefore return any element of \(I_\star\). Sorting costs
\fz{\(\mathcal O(\ell\log(\ell+1))\) comparisons, and the sliding-window scan is linear.
For rational samples and tolerance, the finite-bit runtime additionally depends on the lengths of their encodings.}
\end{proof}

\begin{remark}
The precise choice of computational model is essential for SETH and fine-grained complexity bounds, which concern actual running times rather than merely the distinction between polynomial and exponential complexity.
SETH is usually stated in the standard word-RAM model adopted here.
Our arguments apply similarly to other standard finite-bit models such as deterministic and probabilistic Turing machines and finite-word register or random-access machines, with SETH interpreted in the same model.
Each reduction computes hard networks with parameters within a fixed finite subset of $\Q$ in time polynomial in $\log N$.
We refer to~\cite{compmodel} for a comparison of computation models and to~\cite{sipser13} for further background.
\end{remark}

\subsection{Exponential Time Hypothesis and Strong Exponential Time Hypothesis}\label{sec:seth}
The results of this work rely on the following assumptions from theoretical computer science: 
The \emph{Exponential Time Hypothesis} (ETH) first introduced in~\cite{eth} and the \emph{Strong Exponential Time Hypothesis} (SETH) first introduced in~\cite{seth}.
Both hypotheses give lower bounds on the time complexity for solving SAT problems.
For $k\in\N$, define the exponential-time constant
\[
s_k := \inf\left\{
s>0:
\begin{array}{c}
\text{there exists a deterministic algorithm deciding}\\
\text{$k$-SAT in time }\mathcal O(2^{sn}\operatorname{poly}(n))
\end{array}
\right\}.
\]

\fz{Checking all assignments solves $k$-SAT in time $\order(2^n\operatorname{poly}(n))$, so this infimum is well-defined.}
This also shows $s_k \leq 1$ for all $k \in \N$.
In contrast, $1$SAT and $2$SAT have polynomial time algorithms~\cite{2sat}, implying $s_1 = 0 = s_2$.
\fz{Sch\"oning's randomized local-search algorithm~\cite{walkSat} admits a deterministic derandomization~\cite{MS:FullDerandomizationSchonings2011}.
For every $\varepsilon>0$ and $k\geq3$, its runtime is
\[
    \mathcal O\bigl((2(1-1/k)+\varepsilon)^n\operatorname{poly}(n)\bigr).
\]}
Hence, $s_k\leq\log_2(2(1-1/k))<1$. 

The ETH states that there is no subexponential time algorithm for solving $k$-SAT for any $k \geq 3$: 
\begin{align}\tag{ETH}\label{ETH}
\text{For all $k \geq 3$, we have $s_k > 0$.}
\end{align}
The nestedness $\mathcal{S}_k(n)\subseteq \mathcal{S}_{k+1}(n)$ shows $s_k\leq s_{k+1}$ and hence the limit $\lim_{k \to \infty} s_k\leq 1$ exists.
The Strong Exponential Time Hypothesis (SETH) is the stronger statement
\begin{align}\tag{SETH}\label{SETH}
   \text{The numbers $s_k$ satisfy $\lim_{k \to \infty} s_k = 1$.}
\end{align}
Informally, SETH states that brute-force search is essentially the best possible algorithm for solving $k$-SAT for large $k$.
SETH and, in particular, ETH are widely believed to hold and would imply $P\neq NP$, even though they are still open.
It is known~\cite{IMPAGLIAZZO2001512} that $s_3>0$ follows already if $s_k>0$ for some $k\geq3$.
Hence SETH implies ETH.

This paper focuses on deterministic algorithms, and all reductions used below are deterministic.
However, the same arguments also apply to randomized numerical algorithms: in that case, the deterministic reduction yields a randomized $k$-SAT algorithm with the same success probability.
The corresponding lower bounds therefore hold under the analogous randomized version of SETH, also referred to as randomized SETH.
\fz{Precisely, replace deterministic decision algorithms in the definitions of $s_k$ and $\widetilde s_k$ below by randomized algorithms that succeed with probability at least $2/3$ on every input and have a deterministic worst-case runtime bound.
Randomized SETH asserts that the resulting constants $s_k$ tend to one as $k\to\infty$; whenever a randomized lower bound is invoked, the exponential-time constants are understood in this sense.}

\bigskip
It will be convenient to restrict $k$-SAT to formulas with a number of clauses that is linear in the number of variables.
The sparsification lemma shows that this restriction does not affect SETH.

\begin{lemma}[{Sparsification~\cite[Theorem~1]{IMPAGLIAZZO2001512}}]
\label{lem:sparsification}
	For every $k\in\N$ and $\eps>0$, there exists a constant $C=C(k,\eps)>0$ such that every
	$\alpha\in\mathcal S_k(n)$ admits formulas $\alpha_1,\ldots,\alpha_M\in\mathcal S_k(n)$ with $M\leq 2^{\eps n}$ and at most $Cn$ clauses each such that $\alpha$ is satisfiable if and only if at least one $\alpha_j$, $j\in\CSrange M$, is satisfiable.
	Moreover, the formulas $\alpha_1,\ldots,\alpha_M$ can be computed in time $\mathcal O(\operatorname{poly}(n)2^{\eps n})$.
\end{lemma}

For any $k\geq3$, let $C_k:=C(k,1/k)$ be the constant from \Cref{lem:sparsification} and define
\begin{align*}
    \widetilde{\mathcal S}_k(n)
    :=
    \left\{
        \alpha\in\mathcal S_k(n)
        :
        \#\alpha\leq C_k n
    \right\}.
\end{align*}
We further define the exponential-time constant for sparse $k$-SAT by
\begin{align*}
    \widetilde s_k
    :=
    \inf\left\{
        t>0 :
        \begin{array}{c}
        \text{there exists a deterministic algorithm deciding}
        \\
        \text{satisfiability on $\widetilde{\mathcal S}_k(n)$ in time }
        \mathcal O(2^{tn})
        \end{array}
    \right\}.
\end{align*}

\begin{lemma}\label{lem:sparse_seth}
    For every $k\geq 3$,
    \begin{align*}
        s_k-\frac1k \leq \widetilde s_k \leq s_k.
    \end{align*}
    In particular, SETH is equivalent to $\lim_{k\to\infty} \widetilde s_k=1$.
\end{lemma}
As a consequence, we will frequently refer to the statement $\lim_{k\to\infty} \widetilde s_k=1$ as SETH in the following.
\begin{proof}
	For every $\rho>0$, the definition of $s_k$ yields an algorithm on $\mathcal S_k(n)$ with runtime $\mathcal O(2^{(s_k+\rho/2)n}\operatorname{poly}(n))$.
	Restricting it to $\widetilde{\mathcal S}_k(n)$ and absorbing the polynomial factor into $2^{\rho n/2}$ provides $\widetilde s_k\leq s_k+\rho$.
	Since $\rho>0$ is arbitrary, $\widetilde s_k\leq s_k$ follows such that $\lim_{k\to \infty}\widetilde s_k=1$ implies SETH.

    For the converse direction, let $\eta>0$.
    By the definition of $\widetilde s_k$, there exists an algorithm
    $\mathcal A$ deciding satisfiability of formulas in
    $\widetilde{\mathcal S}_k(n)$ with runtime $\mathcal O\left(2^{(\widetilde s_k+\eta)n}\right)$.

    We now apply the sparsification lemma, Lemma~\ref{lem:sparsification}, with
    $\varepsilon=1/k$.
    Thus, for every $\alpha\in\mathcal S_k(n)$, in runtime $\mathcal O\left(\operatorname{poly}(n)2^{n/k}\right)$ we can construct formulas
    $\alpha_1,\ldots,\alpha_m$ such that $m\leq 2^{n/k}$, each $\alpha_j$ has at most $C_kn$ clauses and hence belongs to
    $\widetilde{\mathcal S}_k(n)$, and $\alpha$ is satisfiable if and
    only if at least one of the formulas $\alpha_j$ is satisfiable.

    We therefore decide satisfiability of $\alpha$ by first applying
    the sparsification lemma and then running $\mathcal A$ on each
    $\alpha_j$. The total runtime is bounded by
    \begin{align*}
        &\mathcal O\left(\operatorname{poly}(n)2^{n/k}\right)
        +
        2^{n/k}
        \mathcal O\left(2^{(\widetilde s_k+\eta)n}\right) =
        \mathcal O\left(
            \operatorname{poly}(n)2^{n/k}
            +
            2^{(\widetilde s_k+1/k+\eta)n}
        \right).
    \end{align*}
    Since $\operatorname{poly}(n)=\mathcal O(2^{\rho n})$ for every
    $\rho>0$, the above
    algorithm has runtime $\mathcal O\left( 2^{(\widetilde s_k+1/k+\eta +\rho)n}\right)$
    for every $\rho>0$.
    By the definition of $s_k$, we have
    \begin{align*}
        s_k
        \leq
        \widetilde s_k+\frac1k+\eta+\rho.
    \end{align*}
    Since $\eta,\rho>0$ are arbitrary, we obtain $s_k- 1/k\leq\widetilde s_k$. Under SETH, $s_k\to1$ as $k\to\infty$, while $1/k\to0$.
    Therefore, $\lim_{k\to\infty}\widetilde s_k=1$.
\end{proof}
\fz{\Cref{lem:sparse_seth} also holds for the randomized exponential-time constants.
For each sparse formula $\alpha_j$, repeat a bounded-error solver independently $R$ times and take the majority answer, where $R$ is the smallest odd integer with $R\geq18\log(3m)$ and $\log$ denotes the natural logarithm.
Each answer is wrong with probability at most $1/3$, so Hoeffding's inequality bounds the probability of a wrong majority by $\exp(-R/18)\leq1/(3m)$.
The union bound over the $m$ formulas therefore gives probability at least $2/3$ that all decisions, and hence their disjunction, are correct.
Since $m\leq2^{n/k}$, we have $R\leq18(\log3+(n/k)\log2)+2=\mathcal O(n)$.
This additional polynomial factor can be absorbed into an arbitrarily small exponential factor, exactly as in the proof above.}

\subsection{Neural networks}
For any subset $\theta\subseteq\R$, we define $\nn(d,w,L,\theta)$ as the set of feed-forward neural networks with $d\in\N$ inputs, maximal width $w\in\N$, maximal depth $L\in\N$, and parameters in $\theta$.
A network $\Phi\in \nn(d,w,L,\theta)$ is a tuple
\begin{align*}
    \Phi
    =
    \left(
        \big(\bW_1,\ldots,\bW_{L'}\big),
        \big(\bb_1,\ldots,\bb_{L'}\big)
    \right)
\end{align*}
with $L'\leq L$, where $\bW_i\in \theta^{w_{i}\times w_{i-1}}$ and $\bb_i\in \theta^{w_{i}}$ for all $i\in\CSrange{L'}$ with $w_0=d$ and $1\leq w_0,\dots,w_{L'}\leq w$.
The network $\Phi$ has $L'-1$ hidden layers and depth $L'$.
Figure~\ref{fig:cnf_to_nn} displays a feed-forward neural network with $d=4$ inputs, three hidden layers ($L' = 4$), and width ten ($w = 10$). All networks in this paper use the ReLU activation $\relu(t):=\max\{0,t\}$, applied componentwise.
The realization $\RR_\Phi:\mathbb{R}^d\to \mathbb{R}^{w_{L'}}$ of $\Phi \in \nn(d,w,L,\theta)$ is the function defined, for any $x\in\mathbb{R}^d$, recursively by
\begin{align*}
	\RR_{\Phi, 0}(x) := x, \qquad \RR_{\Phi, \ell}(x) &:= \relu(\bW_{\ell}\RR_{\Phi, \ell-1}(x)+ \bb_{\ell})\quad\text{for all }\ell\in\CSrange{L'-1},\\
    \RR_\Phi(x) &:= \bW_{L'} \RR_{\Phi, L'-1}(x) + \bb_{L'}.
\end{align*}

\begin{remark}
\label{rem:network_operations}
If the output dimension of a ReLU network \(\Phi\) agrees with the
input dimension of a ReLU network \(\Psi\), their composition is again
a ReLU network and realizes
\(\mathcal R_{\Psi\circ\Phi}=\mathcal R_\Psi\circ\mathcal R_\Phi\).
With the convention on depth used in this paper, the composed depth is
the sum of the two depths minus one. Networks with the same input and
output dimensions can be parallelized and added after padding the
shallower network with identity layers. The resulting depth is the
maximum of the original depths, and the widths and layerwise
sparsities increase at most additively\fz{ up to fixed factors from the identity layers}. Finally,
\(\sigma(t)-\sigma(-t)=t\) provides identity networks of arbitrary
depth. All these constructions are explicit and require only
polynomial time; see also~\cite[Appendix~B]{neuralcalc}.
\end{remark}

\fz{The following three lemmas give an explicit polynomial-time conversion of CNF formulas into neural networks, based on the piecewise affine ramp function in Figure~\ref{fig:ramp}.}

\noindent\begin{minipage}[t]{0.55\textwidth}
\vspace{0pt}
\begin{lemma}\label{lem:r_with_relu}
Let $\eta\in(0,1/2)$ and set $a_\eta := {1}/{(1-2\eta)}$. For $t\in\R$, define
\begin{align*}
    r_\eta(t) :=
    a_\eta
    \Big(
        \sigma(2t-1)
        -
        \sigma(2t-2+2\eta)
    \Big),
\end{align*}
where $\sigma(t)=\max\{0,t\}$. Then $r_\eta:\R\to[0,1]$ satisfies $r_\eta(t) = 0$ for all $t \leq 1/2$ and $r_\eta(t) = 1$ for all $t \geq 1 - \eta$. Moreover, $r_\eta$ is the realization of a ReLU network of
depth $2$ and width $2$.
\end{lemma}
\end{minipage}\hfill
\begin{minipage}[t]{0.42\textwidth}
\vspace{0pt}
\centering
\includegraphics[width=\linewidth]{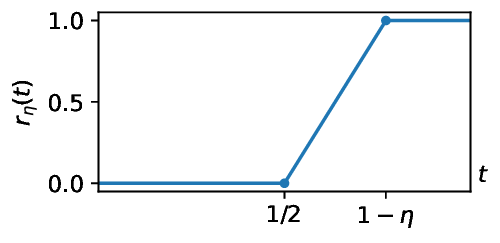}
\setlength{\captionindent}{0pt}
\setlength{\abovecaptionskip}{-20pt}
\captionof{figure}{Ramp function $r_\eta$.}
\label{fig:ramp}
\end{minipage}
\par\medskip
\begin{proof}
If $t\leq1/2$, both arguments of the ReLU functions are
non-positive and hence $r_\eta(t)=0$. If $1/2<t<1-\eta$, only the first ReLU is active, and therefore $r_\eta(t)=a_\eta(2t-1)\in(0,1)$. Finally, if $t\geq 1-\eta$, both ReLU functions are active and
\begin{align*}
    r_\eta(t) =
    a_\eta
    \big(
        (2t-1)-(2t-2+2\eta)
    \big) =
    a_\eta(1-2\eta)
    =1.
\end{align*}
The network realization follows directly from the displayed formula.
\end{proof}

\begin{figure}[tbp]
    \centering
    \begin{tikzpicture}[shorten >=1pt,->,draw=black!50]
        \tikzstyle{every pin edge}=[<-,shorten <=1pt]
        \tikzstyle{neuron}=[circle,fill=black!25,minimum size=17pt,inner sep=0pt]
        \tikzstyle{input neuron}=[neuron, fill=green!50]
        \tikzstyle{output neuron}=[neuron, fill=red!50]
        \tikzstyle{hidden neuron}=[neuron, fill=blue!50]
        \tikzstyle{constant neuron}=[neuron, fill=orange!70]
        \tikzstyle{annot} = [text width=4em, text centered]
        % Parameters
        \def\clauseSep{4.3}
        \def\pairSep{1.3}
        \def\inPairSep{0.5}
        \def\layersep{3.0}

        % Input layer
        \foreach \name / \y in {1,...,4}
            \node[input neuron,font=\scriptsize] (I-\name) at (0,-\y) {$x_{\name}$};

        % Calculate centers
        \def\inputCenter{-2.5}
        \pgfmathsetmacro{\hiddenTop}{0 + 0.5*\inPairSep}
        \pgfmathsetmacro{\hiddenBottom}{-\clauseSep -1*\pairSep -0.5*\inPairSep} % 2 pairs in 2nd clause
        \pgfmathsetmacro{\hiddenCenter}{0.5*(\hiddenTop+\hiddenBottom)}
        \pgfmathsetmacro{\hiddenShift}{\inputCenter-\hiddenCenter}

        % First hidden layer: pairs for each literal in each clause
        % Clause 1: 3 literals (x1, x2, x3)
        \foreach \pair [count=\p from 0] in {1,2,3} {
            \pgfmathsetmacro{\y}{-0*\clauseSep -\p*\pairSep}
            \node[hidden neuron] (H1-1-\pair-1) at (\layersep, {\y + 0.5*\inPairSep + \hiddenShift}) {};
            \node[hidden neuron] (H1-1-\pair-2) at (\layersep, {\y - 0.5*\inPairSep + \hiddenShift}) {};
        }
        % Clause 2: 2 literals (x2, x4)
        \foreach \pair [count=\p from 0] in {1,2} {
            \pgfmathsetmacro{\y}{-1*\clauseSep -\p*\pairSep}
            \node[hidden neuron] (H1-2-\pair-1) at (\layersep, {\y + 0.5*\inPairSep + \hiddenShift}) {};
            \node[hidden neuron] (H1-2-\pair-2) at (\layersep, {\y - 0.5*\inPairSep + \hiddenShift}) {};
        }

        % Second hidden layer: two clause neurons and one constant neuron
        \foreach \clause/\numPairs in {1/3,2/2} {
            \pgfmathsetmacro{\clauseTop}{-(\clause-1)*\clauseSep + 0.5*\inPairSep + \hiddenShift}
            \pgfmathsetmacro{\clauseBottom}{-(\clause-1)*\clauseSep -(\numPairs-1)*\pairSep -0.5*\inPairSep + \hiddenShift}
            \pgfmathsetmacro{\clauseCenter}{0.5*(\clauseTop+\clauseBottom)}
            \node[hidden neuron,font=\scriptsize] (H2-\clause-1) at (2*\layersep, {\clauseCenter + 0.55}) {$u_{\clause}$};
            \node[hidden neuron,font=\scriptsize] (H2-\clause-2) at (2*\layersep, {\clauseCenter}) {$v_{\clause}$};
            \node[constant neuron,font=\scriptsize] (H2-\clause-3) at (2*\layersep, {\clauseCenter - 0.55}) {$e_{\clause}$};
        }

        % Third hidden layer: pair at center of all hidden nodes
        \pgfmathsetmacro{\layerCenter}{\hiddenCenter+\hiddenShift}
        \node[hidden neuron,font=\scriptsize] (H3-1) at (3*\layersep, {\layerCenter + 0.3}) {$p$};
        \node[hidden neuron,font=\scriptsize] (H3-2) at (3*\layersep, {\layerCenter - 0.3}) {$q$};

        % Output layer: single node at center
        \node[output neuron,font=\scriptsize] (O) at (4*\layersep, {\layerCenter}) {$\mathcal R$};

        % Connect input to first hidden layer: each pair to its variable
        % Clause 1: literals x1, x2, x3
        \foreach \pair/\var in {1/1,2/2,3/3}
            \foreach \j in {1,2}
                \draw (I-\var) -- (H1-1-\pair-\j);
        % Clause 2: literals x2, x4
        \foreach \pair/\var in {1/2,2/4}
            \foreach \j in {1,2}
                \draw (I-\var) -- (H1-2-\pair-\j);

        % Connect first hidden to second hidden layer
        \foreach \clause/\numPairs in {1/3,2/2}
            \foreach \pair in {1,...,\numPairs}
                \foreach \j in {1,2}
                    \foreach \k in {1,2}
                        \draw (H1-\clause-\pair-\j) -- (H2-\clause-\k);

        % Connect all u_C, v_C, and e_C neurons to p and q
        \foreach \clause in {1,2}
            \foreach \j in {1,2,3}
                \foreach \k in {1,2}
                    \draw (H2-\clause-\j) -- (H3-\k);

        % Connect third hidden to output
        \foreach \j in {1,2}
            \draw (H3-\j) -- (O);

        \node[draw, rounded corners, fit=(I-1)(I-2)(I-3)(I-4), inner sep=5mm, label=above:{input layer}] {};
        \node[draw, rounded corners, fit=(H1-1-1-1)(H1-1-3-2)(H1-2-1-1)(H1-2-2-2), inner sep=5mm, label=above:{hidden layer 1}] {};
        \node[draw, rounded corners, fit=(H2-1-1)(H2-1-3)(H2-2-1)(H2-2-3), inner sep=5mm, label=above:{hidden layer 2}] {};
        \node[draw, rounded corners, fit=(H3-1)(H3-2), inner sep=5mm, label=above:{hidden layer 3}] {};
        \node[draw, rounded corners, fit=(O), inner sep=5mm, label=above:{output layer}] {};
    \end{tikzpicture}
    \caption{The connectivity graph of $\Psi_{\alpha,\eta}$ from
    Lemma~\ref{lem:CNF_evaluator}, which is the growing-dimensional
    network used in Proposition~\ref{prop:sat_to_nn}, for
    $\alpha=(\gamma_1(x_1)\lor\gamma_2(x_2)\lor\gamma_3(x_3))
    \land(\gamma_4(x_2)\lor\gamma_5(x_4))$. Each literal occurrence
    contributes the pair $h^+,h^-$ in the first hidden layer. For
    each clause $C_j$, the second hidden layer contains $u_j$, $v_j$,
    and the constant neuron $e_j=1$ (orange, with no incoming edge).
    The third hidden layer contains $p,q$, and the output is
    $a_\eta(p-q)$. Thus the layer widths are $10$, $6$, $2$, and $1$,
    exactly as in the lemma. Here
    $\gamma_i\in\{\mathrm{id},\neg\}$.}
    \label{fig:cnf_to_nn}
\end{figure}

\begin{lemma}\label{lem:CNF_evaluator}
Let $k\in\N$, $\eta\in(0,1/2)$, and
$\alpha\in\widetilde{\mathcal S}_k(n)$. Set
\begin{align*}
    m:=\#\alpha,
    \qquad
    \ell(\alpha):=\sum_{C\in\alpha}\#C,
    \qquad
    a_\eta:=\frac{1}{1-2\eta}.
\end{align*}
Then one can construct a ReLU neural network
$\Psi_{\alpha,\eta}$ with input dimension $n$ and realization
\begin{align}\label{eq:robust_CNF_evaluator}
    \mathcal R_{\Psi_{\alpha,\eta}}(x)
    =
    r_\eta\left(
        1-
        \sum_{C\in\alpha}
        \left[
            1-
            r_\eta\left(
                \sum_{(i,\gamma)\in C}
                r_\eta(\gamma(x_i))
            \right)
        \right]
    \right) = r_\eta\left(\left(
        \sum_{C\in\alpha}
            r_\eta\left(
                \sum_{(i,\gamma)\in C}
                r_\eta(\gamma(x_i))\right) - (m - 1)
            \right)
    \right)
\end{align}
The network has depth $4$ and hidden-layer widths
$2\ell(\alpha)$, $3m$ and $2$. In particular,
\begin{align*}
    \operatorname{width}(\Psi_{\alpha,\eta})
    \leq \max\{n,2\ell(\alpha),3m,2\}.
\end{align*}
Moreover, $\mathcal R_{\Psi_{\alpha,\eta}}(z) = \mathcal I_\alpha(z)$ holds for all $z\in\{0,1\}^n$. The numbers of non-zero entries in its weight matrices
$\bW_1,\ldots,\bW_4$ are bounded by
\begin{align*}
    \#\operatorname{supp}(\bW_1)=2\ell(\alpha),\qquad
\#\operatorname{supp}(\bW_2)=4\ell(\alpha),\qquad
\#\operatorname{supp}(\bW_3)=6m,\qquad
\#\operatorname{supp}(\bW_4)=2.
\end{align*}
All weights and biases belong to the finite set
\begin{align}\label{eq:theta_eval}
    \theta_{\rm eval,\eta}
    :=
    \left\{
    0,\pm1,\pm2,\pm a_\eta,\pm2a_\eta,
    \pm2\eta,\pm(2-2\eta)
    \right\}.
\end{align}

The construction can be carried out in polynomial time in the size
of $\alpha$, and $0\leq\mathcal R_{\Psi_{\alpha,\eta}}(x)\leq1$ holds for all $x\in[0,1]^n$. If $\alpha$ is unsatisfiable, then $\mathcal R_{\Psi_{\alpha,\eta}}(x) = 0$ holds for all $x \in [0,1]^n$. If $\alpha$ is satisfiable and
$z\in\{0,1\}^n$ is a satisfying assignment, then $\mathcal R_{\Psi_{\alpha,\eta}}(x)=1$ holds for every $x\in[0,1]^n$ satisfying $|x_i-z_i|\leq\eta$ for all $i\in\CSrange n$. Consequently,
\begin{align*}
    \alpha\text{ satisfiable}
    \quad\Longrightarrow\quad
    \int_{[0,1]^n}
    \mathcal R_{\Psi_{\alpha,\eta}}(x)\,dx
    \geq\eta^n.
\end{align*}
\end{lemma}
\begin{proof}
We construct the network explicitly. By the definition of
$\mathcal S_k(n)$, we have $m\geq1$ and every clause contains at least
one literal. For every occurrence of a literal $(i,\gamma)$ in a clause
$C\in\alpha$, the first hidden layer contains two neurons. If $\gamma=\mathrm{id}$, these neurons compute
\begin{align*}
    h^+_{C,i,\gamma}(x) :=\sigma(2x_i-1) \quad \text{and} \quad 
    h^-_{C,i,\gamma}(x) :=\sigma(2x_i-2+2\eta).
\end{align*}
If $\gamma=\neg$, they compute
\begin{align*}
    h^+_{C,i,\gamma}(x):=\sigma(-2x_i+1) \quad \text{and} \quad
    h^-_{C,i,\gamma}(x):=\sigma(-2x_i+2\eta).
\end{align*}
Thus, in both cases $r_\eta(\gamma(x_i)) = a_\eta \left( h^+_{C,i,\gamma}(x) - h^-_{C,i,\gamma}(x) \right)$. There are exactly $2\ell(\alpha)$ neurons in this layer.
Each neuron is connected to exactly one input variable. Hence
$\bW_1$ has exactly $2\ell(\alpha)$ non-zero entries.

For every clause $C\in\alpha$, the second hidden layer contains two
neurons
\begin{align*}
    u_C(x)
    &:=
    \sigma\left(
        2a_\eta
        \sum_{(i,\gamma)\in C}
        \left(
            h^+_{C,i,\gamma}(x)
            -
            h^-_{C,i,\gamma}(x)
        \right)
        -1
    \right),\\
    v_C(x)
    &:=
    \sigma\left(
        2a_\eta
        \sum_{(i,\gamma)\in C}
        \left(
            h^+_{C,i,\gamma}(x)
            -
            h^-_{C,i,\gamma}(x)
        \right)
        -2+2\eta
    \right).
\end{align*}
By the formula from the literal in the first layer we have
\begin{align*}
    c_C(x) :=
    a_\eta(u_C(x)-v_C(x)) =
    r_\eta\left(
        \sum_{(i,\gamma)\in C}
        r_\eta(\gamma(x_i))
    \right).
\end{align*}
In addition to $u_C$ and $v_C$, the second hidden layer contains for
every clause $C\in\alpha$ a constant neuron
$e_C(x):=\sigma(1)=1$. Hence, this layer has width $3m$. The constant
neurons have no non-zero incoming weights, so
$\#\operatorname{supp}(\bW_2)=4\ell(\alpha)$ remains unchanged.

The third hidden layer consists of two neurons
\[
p(x):=\sigma\left(
1+\sum_{C\in\alpha}
\left(2a_\eta u_C(x)-2a_\eta v_C(x)-2e_C(x)\right)
\right)
\]
and
\[
q(x):=\sigma\left(
2\eta+\sum_{C\in\alpha}
\left(2a_\eta u_C(x)-2a_\eta v_C(x)-2e_C(x)\right)
\right).
\]
Thus, $\#\operatorname{supp}(\bW_3)=6m$. Since
$c_C=a_\eta(u_C-v_C)$ and $e_C=1$, we have
\[
1+\sum_{C\in\alpha}(c_C-e_C)
=
1-\sum_{C\in\alpha}(1-c_C).
\]
Consequently,
\[
a_\eta(p(x)-q(x))
=
r_\eta\left(
1-\sum_{C\in\alpha}(1-c_C(x))
\right).
\]
Taking $a_\eta(p-q)$ as the affine output layer
proves~\eqref{eq:robust_CNF_evaluator}. The output weight matrix has two
non-zero entries.

The displayed formulas also show that every parameter belongs to
$\theta_{\rm eval,\eta}$ from~\eqref{eq:theta_eval}.
In particular, this finite parameter set depends only on $\eta$ and
is independent of $n$, $k$, and $\alpha$. Since $r_\eta$ takes values in $[0,1]$, it follows directly from~\eqref{eq:robust_CNF_evaluator} that $0\leq \mathcal R_{\Psi_{\alpha, \eta}}(x) \leq1$.

We next prove the statement for unsatisfiable formulas.
Let $x\in[0,1]^n$ and choose a closest Boolean vector
$z(x)\in\{0,1\}^n$. Then $|x_i-z_i(x)|\leq 1/2$ holds for every $i\in\CSrange n$. Since $\alpha$ is unsatisfiable, there exists a clause
$C_x\in\alpha$ such that $\gamma(z_i(x))=0$ holds for every $(i,\gamma)\in C_x$. Since both $\mathrm{id}$ and $\neg$ are isometries of $[0,1]$ we have
\begin{align*}
    \gamma(x_i) \leq
    \gamma(z_i(x))
    +
    |\gamma(x_i)-\gamma(z_i(x))| =
    |x_i-z_i(x)|
    \leq\frac12.
\end{align*}
Lemma~\ref{lem:r_with_relu} therefore gives
$r_\eta(\gamma(x_i))=0$ for every $(i,\gamma)\in C_x$. Consequently,
$c_{C_x}(x)=r_\eta(0)=0$. Since $0\leq c_C(x)\leq1$ for every
$C\in\alpha$, we have
\[
    1-\sum_{C\in\alpha}(1-c_C(x))
    \leq
    1-(1-c_{C_x}(x))
    =
    0.
\]
Another application of Lemma~\ref{lem:r_with_relu} therefore gives
$\mathcal R_{\Psi_{\alpha,\eta}}(x)=0$. Since
$x\in[0,1]^n$ was arbitrary, the realization vanishes identically.

Suppose now that $\alpha$ is satisfiable and let
$z\in\{0,1\}^n$ be a satisfying assignment. Let $x\in[0,1]^n$ 
satisfy $|x_i-z_i|\leq\eta$ for all $i\in\CSrange n$. For every clause $C\in\alpha$, there exists a literal
$(i,\gamma)\in C$ with $\gamma(z_i)=1$. Hence, $\gamma(x_i) \geq \gamma(z_i) - |\gamma(x_i)-\gamma(z_i)| \geq 1-\eta$. Lemma~\ref{lem:r_with_relu} gives $r_\eta(\gamma(x_i))=1$, and therefore $c_C(x)=1$. Since this holds
for every clause $C\in\alpha$,
\[
    1-\sum_{C\in\alpha}(1-c_C(x))=1.
\]
Consequently,
$\mathcal R_{\Psi_{\alpha,\eta}}(x)=r_\eta(1)=1$.

The same argument also proves the claimed identity on the Boolean
cube. Let $z\in\{0,1\}^n$. If $\mathcal I_\alpha(z)=0$, then some
clause $C$ is false under $z$, and hence $c_C(z)=0$, which gives
$\mathcal R_{\Psi_{\alpha,\eta}}(z)=0$. If
$\mathcal I_\alpha(z)=1$, then every clause contains a literal which
is true under $z$, and hence $c_C(z)=1$ for every $C\in\alpha$.
Therefore $\mathcal R_{\Psi_{\alpha,\eta}}(z)=1$.

Finally, define
\begin{align*}
    Q_z
    :=
    \left\{
        x\in[0,1]^n:
        |x_i-z_i|\leq\eta
        \text{ for every }i
    \right\}.
\end{align*}
For each coordinate, the corresponding interval has length $\eta$,
and therefore $|Q_z|=\eta^n$. Since the realization is equal to one on $Q_z$,
\begin{align*}
    \int_{[0,1]^n}
    \mathcal R_{\Psi_{\alpha,\eta}}(x)\,dx
    \geq
    |Q_z|
    =
    \eta^n.
\end{align*}
The construction involves only a number of neurons and non-zero
parameters linear in the number of literal occurrences, and hence can
be performed in polynomial time.
\end{proof}

In order to produce neural networks with fixed input dimension, we also introduce the folding map 
\begin{align}\label{eq:folding_map}
    S_\eta(t)
    :=
    \eta^{-1}\Big(
        \sigma(t)-\sigma(t-\eta)
        -\sigma(t-(1-\eta))+\sigma(t-1)
    \Big),
    \qquad t\in[0,1]
\end{align}
for $\eta\in(0,1/2)$.
This piecewise affine function is shown in Figure~\ref{fig:folding_map} and equivalently characterised by
\par\smallskip
\noindent\begin{minipage}[c]{0.55\textwidth}
\begin{align*}
    S_\eta(t)
    =
    \begin{cases}
        t/\eta,&0\leq t\leq\eta,\\
        1,&\eta\leq t\leq1-\eta,\\
        (1-t)/\eta,&1-\eta\leq t\leq1.
    \end{cases}
\end{align*}
\end{minipage}\hfill
\begin{minipage}[c]{0.42\textwidth}
\centering
\includegraphics[width=\linewidth]{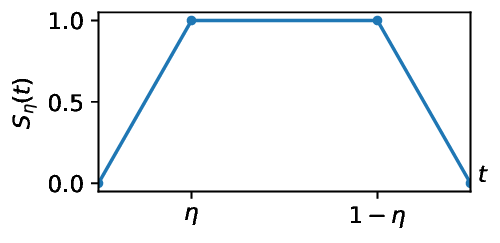}
\setlength{\captionindent}{0pt}
\setlength{\abovecaptionskip}{-20pt}
\captionof{figure}{Folding map $S_\eta$.}
\label{fig:folding_map}
\end{minipage}
\par\medskip
For $q\in\N$, define $f_q:[0,1]\to[0,1]^q$ recursively by
\begin{align}\label{eq:folding_recursion}
    f_1(t):=t,
    \qquad
    f_{q+1}(t)
    :=
    \begin{pmatrix}
        f_q(S_\eta(t))\\
        t
    \end{pmatrix}.
\end{align}

\begin{lemma}\label{lem:fixed_dim_encoder}
Let $\eta\in(0,1/2)$, $s\in\N$, and $n\geq s$. Set $q:=n-s+1$. There exists a ReLU network $\Xi_{n,s,\eta}$ with input dimension
$s$ and output dimension $n$ such that
\begin{align*}
    \mathcal R_{\Xi_{n,s,\eta}}
    (x_1,\ldots,x_s)
    =
    \left(
        x_1,\ldots,x_{s-1},f_q(x_s)
    \right)
\end{align*}
for all $x\in[0,1]^s$. The network can be chosen with depth $n-s+2$, width at most $4n-2s+4$, and at most $16n$ non-zero entries in each weight matrix.
Its final bias vector is zero, and each row of its final weight
matrix contains at most two non-zero entries, all belonging to
$\{-1,1\}$. All weights and biases belong to the finite set
\begin{align*}
    \theta_{\rm CV,\eta}
    :=
    \{0,\pm1,\pm\eta,\pm(1-\eta),\pm\eta^{-1}\}.
\end{align*}

Moreover, for every $z\in\{0,1\}^n$ there exists a measurable set
$Q_{z,s}\subseteq[0,1]^s$ such that $|Q_{z,s}|=\eta^n$ and $\left| \mathcal R_{\Xi_{n,s,\eta}}(x)_i-z_i \right| \leq\eta$ holds for all $x\in Q_{z,s}$ and all $i=1,\ldots,n$.
\end{lemma}
\begin{proof}
We first prove the corresponding property of $f_q$.
For every $q\in\N$ and every $z\in\{0,1\}^q$, there exists an
interval $I_{q,z}\subseteq[0,1]$ of length $\eta^q$ such that $|f_q(t)_i-z_i|\leq\eta$ holds for all $t\in I_{q,z}$ and all $i=1,\ldots,q$. For $q=1$, choose $I_{1,0}=[0,\eta]$ and $I_{1,1}=[1-\eta,1]$.

Suppose the claim holds for $q$ and let
$z=(z',z_{q+1})\in\{0,1\}^{q+1}$.
If $z_{q+1}=0$, define
\begin{align*}
    I_{q+1,z}
    :=
    \{t\in[0,\eta]:S_\eta(t)\in I_{q,z'}\}.
\end{align*}
On $[0,\eta]$, we have $S_\eta(t)=t/\eta$, and therefore, $|I_{q+1,z}|=\eta|I_{q,z'}|=\eta^{q+1}$. Moreover, $|t-z_{q+1}|\leq\eta$, while the induction hypothesis
applied to $S_\eta(t)$ yields closeness to the corner for the first $q$ coordinates. If $z_{q+1}=1$, the same argument applies on
$[1-\eta,1]$, where $S_\eta(t)={(1-t)}/{\eta}$. This proves the claim for $f_{q+1}$.

Now let $z=(z_1,\ldots,z_n)\in\{0,1\}^n$.
For $i=1,\ldots,s-1$, let
\begin{align*}
    J_{z_i}
    :=
    \begin{cases}
        [0,\eta],&z_i=0,\\
        [1-\eta,1],&z_i=1.
    \end{cases}
\end{align*}
Let $I\subseteq[0,1]$ be the interval of length $\eta^q$
corresponding to the Boolean vector $(z_s,\ldots,z_n)$ in the
previous construction, and define $Q_{z,s} := J_{z_1}\times\cdots\times J_{z_{s-1}}\times I$.
Then $|Q_{z,s}| = \eta^{s-1}\eta^q = \eta^n$, because $q=n-s+1$. The required coordinatewise estimate follows
immediately from the definition of $\Xi_{n,s,\eta}$. It remains to prove the network bounds.

The function $S_\eta$ is realized on $[0,1]$ by the depth-$2$,
width-$4$ network in~\eqref{eq:folding_map} with zero output bias. Using composition, parallelization, and identity networks, we construct
a network $\Gamma_q$ realizing $f_q$ recursively from~\eqref{eq:folding_recursion}. An induction shows that $\Gamma_q$ can be chosen with $\operatorname{depth}(\Gamma_q)=q+1$ and $\operatorname{width}(\Gamma_q)\leq4q+2$,
and with at most $16q$ non-zero entries in every weight matrix.

For $q=1$, recall that $f_1(t)=t$. We may take $\Gamma_1$ to be the
depth-$2$ ReLU network $t\mapsto \sigma(t)-\sigma(-t)=t$. It has width $2$, at most two non-zero entries in each weight matrix,
zero final bias, entries in $\{-1,0,1\}$ in its first weight matrix,
and its final weight matrix has two non-zero entries in
$\{-1,1\}$. Thus all asserted bounds hold for $q=1$.

Suppose these properties hold for $\Gamma_q$.
The composition $\Gamma_q\circ S_\eta$ has depth $q+2$.
Parallelizing it with an identity network of depth $q+2$ gives
$\Gamma_{q+1}$ and yields
\begin{align*}
    \operatorname{width}(\Gamma_{q+1}) \leq
    \max\{4q+2,4\}+2 \leq4(q+1)+2.
\end{align*}
All unchanged weight matrices contain at most $16q$ non-zero entries.
The only potentially larger matrix is the matrix obtained by merging
the output layer of the network realizing $S_\eta$ with the first
layer of $\Gamma_q$. Hence, this product has at
most $4(4q+2)=16q+8$ non-zero entries. Parallelization with the identity branch adds at
most four further non-zero entries, and thus, $16q+12\leq16(q+1)$. This proves the induction.

The construction can furthermore be performed so that the first
weight matrix has entries in $\{-1,0,1\}$, the final bias vanishes,
and every row of the final weight matrix has at most two non-zero
entries in $\{-1,1\}$. These properties are preserved by the same
induction. Since the output bias of the network realizing $S_\eta$
is zero, composition creates no new bias values. Thus, all parameters
remain in $\theta_{\rm CV,\eta}$.

Finally, parallelize $\Gamma_q$ with $s-1$ identity networks of depth
$q+1$, one for each of the first $s-1$ input coordinates.
This produces $\Xi_{n,s,\eta}$.
Its width is bounded by $(4q+2)+2(s-1) = 4n-2s+4$, and each weight matrix has at most $16q+4(s-1) \leq16n$
non-zero entries.
The asserted depth is $q+1=n-s+2$. This completes the proof.
\end{proof}

\begin{proposition}\label{prop:sat_to_nn}
Let $k\geq3$, let $\delta\in\Q_{>0}$, and set
$\eta:=(2+\delta)^{-1}$ and $a_\eta:=(1-2\eta)^{-1}$. Define
\[
\theta_{{\rm CNF},\delta}
:=
\big\{
0,\pm1,\pm2,
\pm\eta,\pm(1-\eta),
\pm2\eta,\pm(2-2\eta),
\pm\eta^{-1},
\pm a_\eta,\pm2a_\eta
\big\}.
\]
Furthermore, set
$C_k^{\rm width}:=\max\{6,2kC_k+2\}$ and
$C_k^{\rm sparse}:=\max\{16,4kC_k+2\}$. For every $n\in\N$ and every
$\alpha\in\widetilde{\mathcal S}_k(n)$ there exists a ReLU network
$\Psi_{\alpha,\delta}$ with input dimension $n$ and scalar output
such that
$0\leq\mathcal R_{\Psi_{\alpha,\delta}}\leq1$ on $[0,1]^n$,
its depth is at most $4$, its width is bounded by
$(2kC_k+2)n$, and every weight matrix has at most
$(4kC_k+2)n$ non-zero entries. Moreover, $\mathcal R_{\Psi_{\alpha,\delta}}(z) = \mathcal I_\alpha(z)$ holds for all $z\in\{0,1\}^n$.

For every fixed $s\in\N$ and every $n\geq s$, there exists a ReLU
network $\Phi_{\alpha,s,\delta}$ with input dimension $s$ and scalar
output such that
$0\leq\mathcal R_{\Phi_{\alpha,s,\delta}}\leq1$ on $[0,1]^s$,
its depth is at most $n-s+5$, its width is bounded by
$C_k^{\rm width}n$, and every weight matrix has at most
$C_k^{\rm sparse}n$ non-zero entries. Both networks can be constructed from $\alpha$ in polynomial time,
and all their weights and biases belong to
$\theta_{{\rm CNF},\delta}$.

If $\alpha$ is unsatisfiable, then
$\mathcal R_{\Psi_{\alpha,\delta}}=0$ on $[0,1]^n$ and
$\mathcal R_{\Phi_{\alpha,s,\delta}}=0$ on $[0,1]^s$. If $\alpha$
is satisfiable, then
\[
\int_{[0,1]^n}
\mathcal R_{\Psi_{\alpha,\delta}}(x)\,dx
\geq
(2+\delta)^{-n}
\qquad\text{and}\qquad
\int_{[0,1]^s}
\mathcal R_{\Phi_{\alpha,s,\delta}}(x)\,dx
\geq
(2+\delta)^{-n}.
\]
\end{proposition}

\begin{proof}
For the $n$-dimensional network, apply
Lemma~\ref{lem:CNF_evaluator} with $\eta=(2+\delta)^{-1}$.
In the notation of that lemma, set
$\Psi_{\alpha,\delta}:=\Psi_{\alpha,\eta}$; the subscript $\delta$
emphasizes the parameter used in the complexity statements below.
Since $\alpha\in\widetilde{\mathcal S}_k(n)$, we have
$\#\alpha\leq C_k n$ and
$\ell(\alpha)\leq k C_k n$. Hence,
\[
\operatorname{width}(\Psi_{\alpha,\delta})
\leq
\max\{n,2k C_k n,3C_k n,2\}
\leq
(2kC_k+2)n,
\]
where we used $k\geq3$ and $n\geq1$. Similarly,
Lemma~\ref{lem:CNF_evaluator} gives
\[
\max_{\ell}
\#\operatorname{supp}(\bW_\ell)
\leq
\max\{2k C_k n,4k C_k n,6C_k n,2\}
\leq
(4kC_k+2)n.
\]
The depth bound, the range
$0\leq\mathcal R_{\Psi_{\alpha,\delta}}\leq1$, and the identity
$\mathcal R_{\Psi_{\alpha,\delta}}(z)=\mathcal I_\alpha(z)$ on
$\{0,1\}^n$ follow directly from
Lemma~\ref{lem:CNF_evaluator}.

Now fix $s\leq n$. Let $\Xi_{n,s,\eta}$ be the network from
Lemma~\ref{lem:fixed_dim_encoder} and define $\Phi_{\alpha,s,\delta} := \Psi_{\alpha,\delta}\circ\Xi_{n,s,\eta}$. Since $\Xi_{n,s,\eta}$ has depth $n-s+2$ and
$\Psi_{\alpha,\delta}$ has depth $4$, the standard composition rule
gives
\[
\operatorname{depth}(\Phi_{\alpha,s,\delta})
\leq
(n-s+2)+4-1
=
n-s+5.
\]
The encoder has width at most $4n-2s+4\leq6n$, while
$\Psi_{\alpha,\delta}$ has width at most $(2kC_k+2)n$.
Consequently, $\operatorname{width}(\Phi_{\alpha,s,\delta}) \leq C_k^{\rm width}n$.

All weight matrices inherited from the encoder contain at most $16n$
non-zero entries. The weight matrices of the CNF evaluator satisfy
the bounds above. It remains to estimate the matrix created by
merging the final affine layer of $\Xi_{n,s,\eta}$ with the first
affine layer of $\Psi_{\alpha,\delta}$. Every row of the first weight matrix of
$\Psi_{\alpha,\delta}$ contains exactly one non-zero entry, and there
are $2\ell(\alpha)$ such rows. By
Lemma~\ref{lem:fixed_dim_encoder}, every row of the final weight
matrix of $\Xi_{n,s,\eta}$ contains at most two non-zero entries.
Therefore, every row of the composed matrix contains at most two
non-zero entries and the matrix contains at most $4\ell(\alpha) \leq 4k C_k n$ non-zero entries. Hence, every weight matrix of
$\Phi_{\alpha,s,\delta}$ contains at most
$C_k^{\rm sparse}n$ non-zero entries.

The parameters of $\Psi_{\alpha,\delta}$ belong to
$\theta_{\rm eval,\eta}$ and those of $\Xi_{n,s,\eta}$ belong to
$\theta_{\rm CV,\eta}$. Moreover, the final bias of
$\Xi_{n,s,\eta}$ is zero and every non-zero entry of its final weight
matrix belongs to $\{-1,1\}$. Thus, merging the final affine layer of
$\Xi_{n,s,\eta}$ with the first affine layer of
$\Psi_{\alpha,\delta}$ creates no new parameter values. Consequently,
all weights and biases of both networks belong to
$\theta_{{\rm CNF},\delta}$.

Both constituent networks can be constructed in polynomial time, and
their composition can therefore also be constructed in polynomial
time. If $\alpha$ is unsatisfiable,
Lemma~\ref{lem:CNF_evaluator} gives
$\mathcal R_{\Psi_{\alpha,\delta}}=0$ on $[0,1]^n$.
Since
$\mathcal R_{\Xi_{n,s,\eta}}([0,1]^s)\subseteq[0,1]^n$, it follows
that
$\mathcal R_{\Phi_{\alpha,s,\delta}}=0$ on $[0,1]^s$.

Suppose that $\alpha$ is satisfiable and let
$z\in\{0,1\}^n$ be a satisfying assignment. By
Lemma~\ref{lem:fixed_dim_encoder}, there exists
$Q_{z,s}\subseteq[0,1]^s$ with $|Q_{z,s}|=\eta^n$ such that $\left| \mathcal R_{\Xi_{n,s,\eta}}(x)_i-z_i \right| \leq\eta$ holds for every $x\in Q_{z,s}$ and every $i \in \CSrange n$. Lemma~\ref{lem:CNF_evaluator} therefore implies
$\mathcal R_{\Phi_{\alpha,s,\delta}}(x)=1$ for every
$x\in Q_{z,s}$. Hence
\[
\int_{[0,1]^s}
\mathcal R_{\Phi_{\alpha,s,\delta}}(x)\,dx
\geq
|Q_{z,s}|
=
\eta^n
=
(2+\delta)^{-n}.
\]
The corresponding $n$-dimensional estimate follows directly from
Lemma~\ref{lem:CNF_evaluator}.
\end{proof}

\FloatBarrier
\section{Quadrature of neural networks}\label{sec:quad}

Proposition~\ref{prop:sat_to_nn} reduces sparse $k$-SAT to the
following gap problem for neural network integrands:
\begin{align*}
    \alpha\text{ unsatisfiable} \Longrightarrow\quad
    \int \mathcal R_{\Phi_\alpha}=0, \qquad
    \alpha\text{ satisfiable} \Longrightarrow\quad
    \int \mathcal R_{\Phi_\alpha}\geq(2+\delta)^{-n}.
\end{align*}
We use this gap first for networks whose input dimension is allowed
to grow and then for networks of fixed input dimension.

\subsection{Quadrature in growing dimension}

\begin{theorem}\label{thm:growing_dim_nn_quadrature_is_hard}
For every $\eps>0$ there exist $C_\eps>0$ and a finite set
$\theta_\eps\subset\Q$ such that the following holds.

Let $C>0$, and suppose that $\cA$ is a deterministic algorithm which,
for every $N\geq2$, accepts every ReLU network $\Phi$ satisfying
\begin{itemize}
    \item the input dimension is at most
    $C_\eps(1+\log_2 N)$;
    \item the width is at most
    $C_\eps(1+\log_2 N)$;
    \item the depth is at most $4$;
    \item every weight matrix contains at most
    $C_\eps(1+\log_2 N)$ non-zero entries;
    \item every weight and bias belongs to $\theta_\eps$,
\end{itemize}
and returns $\cA(\Phi,N)\in\R$ in runtime at most $N$. Suppose, for some $\gamma>0$ and every admissible $(\Phi,N)$, that
\begin{align*}
    \left|
        \cA(\Phi,N)
        -
        \int_{[0,1]^{d_\Phi}}\cR_\Phi(x)\,dx
    \right|
    \leq
    C\norm{\cR_\Phi}{L^2([0,1]^{d_\Phi})}N^{-\gamma}
    +
    C^2N^{-2\gamma},
\end{align*}
where $d_\Phi$ denotes the input dimension of $\Phi$.
Then SETH implies $\gamma\leq 1/2 + \eps$.
\end{theorem}
\begin{proof}
Fix $\eps>0$ and suppose, for contradiction, that $\gamma> 1/2+\eps$. Choose a rational number $\omega$ such that
\begin{align*}
    \max\left\{
        \frac12,\frac{1}{1+2\eps}
    \right\}
    <\omega<1.
\end{align*}
By Lemma~\ref{lem:sparse_seth}, SETH implies
$\widetilde s_k\to1$. Hence, there exists $k\geq3$ such that $\omega<\widetilde s_k$. Choose a rational $\delta>0$ sufficiently small such that $\omega\log_2(2+\delta)<\widetilde s_k$.
We may therefore choose a rational $B>1$ with $(2+\delta)^\omega < B < 2^{\widetilde s_k}$. Set $\theta_\eps:=\theta_{{\rm CNF},\delta}$ with $\theta_{{\rm CNF},\delta}$ from
Proposition~\ref{prop:sat_to_nn}.

For $n\in\N$, define $q_n:=(2+\delta)^{-n}$ and $N_n:=\left\lceil B^n\right\rceil$. Since $\omega>1/2$ we have $B>\sqrt2$. Thus,
\begin{align*}
    n
    \leq
    \frac{\log_2N_n}{\log_2B}
    <
    2\log_2N_n.
\end{align*}

Let $\alpha\in\widetilde{\mathcal S}_k(n)$ and let
$\Psi_{\alpha,\delta}$ be the $n$-dimensional network from
Proposition~\ref{prop:sat_to_nn}. It has depth $4$ and satisfies $\operatorname{width}(\Psi_{\alpha,\delta}) \leq (2kC_k+2)n$ and $\#\operatorname{supp}(\bW_\ell) \leq(4kC_k+2)n$. This shows that,
after choosing $C_\eps$ sufficiently large,
$\Psi_{\alpha,\delta}$ is admissible for
$\cA(\,\cdot\,,N_n)$.

Set
\begin{align*}
    I_\alpha
    :=
    \int_{[0,1]^n}
    \cR_{\Psi_{\alpha,\delta}}(x)\,dx.
\end{align*}
By Proposition~\ref{prop:sat_to_nn}, $I_\alpha=0$ holds for all $\alpha$ that are unsatisfiable and $I_\alpha \geq q_n$ holds for all $\alpha$ that are satisfiable. Moreover, since $0\leq\cR_{\Psi_{\alpha,\delta}}\leq1$, we have $\norm{\cR_{\Psi_{\alpha,\delta}}}{L^2([0,1]^n)}^2 \leq I_\alpha$. Since $\gamma>1/2+\eps$ and
$\omega>1/(1+2\eps)$, we have $\gamma\omega > 1/2$.
Hence, $\rho:=(2+\delta)^{1/2}B^{-\gamma}<1$, and therefore
\begin{align*}
    CN_n^{-\gamma}
    \leq
    CB^{-\gamma n}
    =
    C\rho^n q_n^{1/2}.
\end{align*}
Choose $n_0$ sufficiently large that $C\rho^n\leq 1/8$ holds for all $n\geq n_0$.

If $\alpha$ is unsatisfiable, then $|\cA(\Psi_{\alpha,\delta},N_n)| \leq C^2N_n^{-2\gamma} \leq {q_n}/{64}$. If $\alpha$ is satisfiable, then
\begin{align*}
    C
    \norm{\cR_{\Psi_{\alpha,\delta}}}{L^2}
    N_n^{-\gamma} \leq
    \frac18I_\alpha, \quad \text{and} \quad C^2N_n^{-2\gamma} \leq
    \frac{q_n}{64}
    \leq
    \frac{I_\alpha}{64},
\end{align*}
hold, and hence, $\cA(\Psi_{\alpha,\delta},N_n) \geq {55I_\alpha}/{64} \geq {55q_n}/{64}$. Thus, for $n\geq n_0$, we can decide satisfiability of $\alpha\in\widetilde{\mathcal S}_k(n)$ by returning TRUE if and only if $ \cA(\Psi_{\alpha,\delta},N_n) \geq 2^{-m}$, where $m\in\mathbb{N}$ is such that the threshold satisfies $2^{-m}\in({q_n}/{64}, {55q_n}/{64})$. For $n<n_0$, we use brute force.

The construction of $\Psi_{\alpha,\delta}$ from Proposition~\ref{prop:sat_to_nn} costs $\mathcal{O}(\operatorname{poly}(n))$ time.
Since $m=\mathcal O(n)$ and writing the output is included in the runtime of $\cA$, 
the threshold test using binary shift and integer comparison requires at most $\mathcal O(N_n+n)$ additional work.
Hence the resulting algorithm for sparse $k$-SAT has runtime
\begin{align*}
    \mathcal O(\operatorname{poly}(n)+N_n)
    =
    \mathcal O(B^n)
    =
    \mathcal O\left(2^{(\log_2B)n}\right).
\end{align*}
Since $\log_2B<\widetilde s_k$, this contradicts the definition of
$\widetilde s_k$. Therefore, $\gamma\leq 1/2+\eps$.
\end{proof}

\subsection{Quadrature in fixed dimension}

\begin{theorem}\label{thm:fixed_dim_nn_quadrature_is_hard}
For every $s\in\N$ there exists $C_s>0$ such that for every
$\eps>0$ there exist $C_\eps>0$ and a finite set
$\theta_\eps\subset\Q$ such that the following holds. Let $C>0$, and suppose $\cA$ is a deterministic algorithm which,
for every $N\geq2$, accepts every ReLU network $\Phi$ satisfying
\begin{itemize}
    \item the input dimension is $s$;
    \item the width is at most
    $C_\eps(1+\log_2 N)$;
    \item the depth is at most
    $C_s(1+\log_2 N)$;
    \item every weight matrix contains at most
    $C_\eps(1+\log_2 N)$ non-zero entries;
    \item every weight and bias belongs to $\theta_\eps$,
\end{itemize}
and returns $\cA(\Phi,N)\in\R$ in runtime at most $N$. Suppose, for every admissible $(\Phi,N)$, that
\begin{align*}
    \left|
        \cA(\Phi,N)
        -
        \int_{[0,1]^s}\cR_\Phi(x)\,dx
    \right|
    \leq
    C\norm{\cR_\Phi}{L^2([0,1]^s)}N^{-\gamma}
    +
    C^2N^{-2\gamma}
\end{align*}
holds for some $\gamma > 0$. Then SETH implies $\gamma\leq 1/2+\eps$.
\end{theorem}
\begin{proof}
We follow the proof of
Theorem~\ref{thm:growing_dim_nn_quadrature_is_hard}
with the same choices of
$\omega$, $k$, $\delta$, $B$, $q_n$, and $N_n$.

For
$\alpha\in\widetilde{\mathcal S}_k(n)$,
we use the fixed-dimensional network
$\Phi_{\alpha,s,\delta}$ from
Proposition~\ref{prop:sat_to_nn} instead of
$\Psi_{\alpha,\delta}$.
For $n\geq s$, Proposition~\ref{prop:sat_to_nn} gives
\begin{align*}
    \operatorname{width}(\Phi_{\alpha,s,\delta}) \leq C_k^{\rm width}n, \qquad
    \operatorname{depth}(\Phi_{\alpha,s,\delta})\leq n-s+5, \qquad
    \#\operatorname{supp}(\bW_\ell) \leq C_k^{\rm sparse}n.
\end{align*}
Since $n<2\log_2N_n$, we may choose $C_\eps$ and $C_s$
sufficiently large that these networks are admissible inputs for
$\cA(\,\cdot\,,N_n)$. Proposition~\ref{prop:sat_to_nn} gives exactly the same integral gap,
\begin{align*}
    \alpha\text{ unsatisfiable} \Longrightarrow
    \int_{[0,1]^s}\cR_{\Phi_{\alpha,s,\delta}}=0, \qquad 
    \alpha\text{ satisfiable} \Longrightarrow
    \int_{[0,1]^s}\cR_{\Phi_{\alpha,s,\delta}}
    \geq(2+\delta)^{-n},
\end{align*}
and $0\leq\cR_{\Phi_{\alpha,s,\delta}}\leq1$.
Hence all remaining estimates in the proof of
Theorem~\ref{thm:growing_dim_nn_quadrature_is_hard}
apply verbatim and produce a sparse $k$-SAT algorithm of runtime
$\mathcal O(B^n)$ with $\log_2B<\widetilde s_k$.
This contradicts the definition of $\widetilde s_k$.
\end{proof}

\subsection{Quadrature on other domains in fixed dimension}

We restrict the extension beyond the cube to fixed input dimension.

\begin{lemma}\label{lem:domain_clipping}
Let $s\in\N$, $a\in\Q^s$, and $r\in\Q_{>0}$. Define
$P_{a,r}\colon\R^s\to[0,1]^s$ componentwise by
\begin{align*}
    P_{a,r}(x)_j
    :=
    \sigma\left(\frac{x_j-a_j}{r}\right)
    -
    \sigma\left(\frac{x_j-a_j}{r}-1\right).
\end{align*}
Then $P_{a,r}(x)= {(x-a)}/{r}$ holds for all $x\in a+r[0,1]^s$. Moreover, $P_{a,r}$ is the realization of a ReLU network of depth
$2$, width $2s$, and with at most $2s$ non-zero entries in each
weight matrix. All its parameters belong to a finite subset of
$\Q$ depending only on $a$ and $r$.
\end{lemma}
\begin{proof}
For $t\in\R$, let $p(t):=\sigma(t)-\sigma(t-1)$.
Then
\begin{align*}
    p(t)
    =
    \begin{cases}
        0,&t\leq0,\\
        t,&0\leq t\leq1,\\
        1,&t\geq1.
    \end{cases}
\end{align*}
Hence $p(\R)\subseteq[0,1]$ and
$p(t)=t$ for $t\in[0,1]$. Applying this componentwise to
$t=(x_j-a_j)/r$ proves both claims concerning $P_{a,r}$. The network realization follows directly from the displayed
formula. Each of the $2s$ hidden neurons depends on one input
coordinate, and each output coordinate is the difference of two
hidden neurons.
\end{proof}

\begin{proposition}\label{prop:sat_to_nn_general_domain}
Let $s,k\in\N$ with $k\geq3$, let
$\Omega\subset\R^s$ be a bounded domain, and suppose that there exist
$a\in\Q^s$ and $r\in\Q_{>0}$ such that $a+r[0,1]^s\subseteq\Omega$. Let $\delta\in\Q_{>0}$.

Then there exists a finite set
$\theta_{\Omega,\delta}\subset\Q$ such that, for every
$n\geq s$ and every
$\alpha\in\widetilde{\mathcal S}_k(n)$, one can construct in
polynomial time a ReLU network
$\widehat\Phi_{\alpha,\Omega,\delta}$ of input dimension $s$
satisfying $0\leq \mathcal R_{\widehat\Phi_{\alpha,\Omega,\delta}} \leq 1$ on $\R^s$. Its width and the number of non-zero entries in every weight matrix are bounded by $C_{k,\Omega}n$, and its depth is bounded by $n-s+6$, where $C_{k,\Omega}>0$ is independent of $n$ and $\alpha$.
All weights and biases belong to $\theta_{\Omega,\delta}$.

If $\alpha$ is unsatisfiable, then $\mathcal R_{\widehat\Phi_{\alpha,\Omega,\delta}} \equiv 0$ on $\Omega$. If $\alpha$ is satisfiable, then
\begin{align*}
    \int_\Omega
    \mathcal R_{\widehat\Phi_{\alpha,\Omega,\delta}}(x)\,dx
    \geq
    r^s(2+\delta)^{-n}.
\end{align*}
\end{proposition}
\begin{proof}
Let $P_{a,r}$ be the map from
Lemma~\ref{lem:domain_clipping}, and define $\widehat\Phi_{\alpha,\Omega,\delta} := \Phi_{\alpha,s,\delta}\circ P_{a,r}$, where $\Phi_{\alpha,s,\delta}$ is the fixed-dimensional network from Proposition~\ref{prop:sat_to_nn}.

Since $P_{a,r}(\R^s)\subseteq[0,1]^s$ and
$0\leq\mathcal R_{\Phi_{\alpha,s,\delta}}\leq1$, we have $0\leq \mathcal R_{\widehat\Phi_{\alpha,\Omega,\delta}} \leq 1$. If $\alpha$ is unsatisfiable, Proposition~\ref{prop:sat_to_nn} yields
$\mathcal R_{\Phi_{\alpha,s,\delta}}\equiv0$ on $[0,1]^s$. Hence, $\mathcal R_{\widehat\Phi_{\alpha,\Omega,\delta}} \equiv 0$ on all of $\R^s$.

Suppose $\alpha$ is satisfiable. By
Proposition~\ref{prop:sat_to_nn}, there exists a measurable set
$Q_{\alpha}\subseteq[0,1]^s$ with $|Q_{\alpha}|\geq(2+\delta)^{-n}$ such that $\mathcal R_{\Phi_{\alpha,s,\delta}}=1$ holds on $Q_{\alpha}$. Since $a+rQ_\alpha \subseteq a+r[0,1]^s \subseteq\Omega$, and $P_{a,r}(x)={(x-a)}/{r}$ on $a+r[0,1]^s$, we obtain $\mathcal R_{\widehat\Phi_{\alpha,\Omega,\delta}}(x)=1$ for $x\in a+rQ_\alpha$. Consequently,
\begin{align*}
    \int_\Omega
    \mathcal R_{\widehat\Phi_{\alpha,\Omega,\delta}}(x)\,dx \geq
    |a+rQ_\alpha| =
    r^s|Q_\alpha| \geq
    r^s(2+\delta)^{-n}.
\end{align*}

The network bounds follow from
Proposition~\ref{prop:sat_to_nn},
Lemma~\ref{lem:domain_clipping}, and the standard composition rule.
In particular, the clipping network has depth $2$, so $\operatorname{depth} (\widehat\Phi_{\alpha,\Omega,\delta}) \leq 2+(n-s+5)-1 = n-s+6$. All additional parameters depend only on the fixed rational numbers
$a$ and $r$, which proves the existence of
$\theta_{\Omega,\delta}$.
\end{proof}

\begin{theorem}\label{thm:general_domain_quadrature}
Let $s\in\N$ and let $\Omega\subset\R^s$ be a bounded domain containing a cube
$a+r[0,1]^s$ for some $a\in\Q^s$ and $r\in\Q_{>0}$.
There exists $C_{s,\Omega}>0$ such that for every $\eps>0$
there exist $C_{\eps,\Omega}>0$ and a finite set
$\theta_{\eps,\Omega}\subset\Q$ such that the following holds.

Suppose $C>0$ and $\cA$ is a deterministic algorithm which accepts
the same class of networks as in
Theorem~\ref{thm:fixed_dim_nn_quadrature_is_hard}, with $[0,1]^s$
replaced by $\Omega$\fz{ and the constants $C_s$, $C_\eps$, and $\theta_\eps$ replaced by $C_{s,\Omega}$, $C_{\eps,\Omega}$, and $\theta_{\eps,\Omega}$, respectively}, and has runtime at most $N$.

If, for some $\gamma>0$ we have
\begin{align*}
    \left|
        \cA(\Phi,N)
        -
        \int_\Omega\cR_\Phi(x)\,dx
    \right|
    \leq
    C\norm{\cR_\Phi}{L^2(\Omega)}N^{-\gamma}
    +
    C^2N^{-2\gamma}
\end{align*}
for every admissible $(\Phi,N)$, then SETH implies $\gamma\leq 1/2+\eps$.
\end{theorem}
\begin{proof}
The proof of
Theorem~\ref{thm:fixed_dim_nn_quadrature_is_hard}
applies with
$\widehat\Phi_{\alpha,\Omega,\delta}$ from
Proposition~\ref{prop:sat_to_nn_general_domain}.
The only change is that the satisfiable case now gives the gap $q_n := r^s(2+\delta)^{-n}$. The prefactor $r^s>0$ is independent of $n$ and therefore does not
affect any exponential rate in the proof.
Moreover,
\begin{align*}
    \norm{
        \mathcal R_{\widehat\Phi_{\alpha,\Omega,\delta}}
    }{L^2(\Omega)}^2
    \leq
    \int_\Omega
    \mathcal R_{\widehat\Phi_{\alpha,\Omega,\delta}} dx,
\end{align*}
since the realization takes values in $[0,1]$.
All remaining estimates are unchanged.
\end{proof}

\subsection{Algebraic optimality of random-bit Monte Carlo}
\label{sec:random_bit_mc}

Random-bit quadrature and its information complexity are studied
in~\cite{HNP:HowManyRandom2004}.
For our network classes, the classical Monte Carlo rate also holds,
up to logarithmic factors, in the total word-RAM runtime.
We give the short additional argument accounting for exact network
evaluation and finite-bit sampling. 

\begin{theorem}
\label{thm:random_bit_mc}
Fix a finite nonempty set $\theta\subset\Q$ and $\kappa\geq1$.
There exist constants $C>0$, $r,N_0\in\N$ with $N_0\geq2$, and a uniform random-bit algorithm \fz{$\cA_{\rm MC}$ in the sense of \Cref{sec:algorithms}} with the following property.
For every runtime budget $N\in\N$ with $N\geq N_0$ and every scalar-valued ReLU network $\Phi$ with parameters in $\theta$, input dimension $d$, width $w$, and depth $L$ satisfying
\begin{align*}
    \max\{d,w,L\}\leq\kappa(1+\log_2N),
\end{align*}
the algorithm returns $\cA_{\rm MC}(\Phi,N)\in\Q$ with word-RAM runtime at most $N$ for every outcome of its random bits.
For each such input, with probability at least $2/3$, it satisfies
\begin{align}\label{eq:random_bit_mc_error}
    \left|\cA_{\rm MC}(\Phi,N)-\int_{[0,1]^d}\realization_\Phi(x)\,dx\right|
    \leq
    C\norm{\realization_\Phi}{L^2([0,1]^d)}
    \frac{(1+\log_2N)^{r/2}}{\sqrt N}
    +\frac{C}{N}.
\end{align}
Consequently, for every exponent $0<\gamma<1/2$, there exists $C_\gamma>0$ such that, for every such input, with probability at least $2/3$,
\begin{align}\label{eq:random_bit_mc_mixed}
    \left|\cA_{\rm MC}(\Phi,N)-\int_{[0,1]^d}\realization_\Phi(x)\,dx\right|
    \leq
    C_\gamma\norm{\realization_\Phi}{L^2([0,1]^d)}N^{-\gamma}
    +C_\gamma^2N^{-2\gamma}.
\end{align}
The constants $C$, $r$, and $N_0$ depend only on $\theta$ and $\kappa$, while $C_\gamma$ may also depend on $\gamma$.
\end{theorem}

\begin{proof}
Write $f:=\realization_\Phi$ and $I(f):=\int_{[0,1]^d}f(x)\,dx$.
Fix integers $b,q\geq1$ such that $|a|\leq b$ and $qa\in\Z$
for every $a\in\theta$.
\bg{Throughout the proof, the $\mathcal O$-notation hides constants that depend exclusively on $\theta$ and $\kappa$.}

Since each affine layer has at most $w$ inputs (recall $d\leq w$) and ReLU is
$1$-Lipschitz, both the maximum norm of $f$ and its Lipschitz constant
in the maximum norm are bounded by $K:=[b(w+1)]^L$.
Choose
\[
    p:=\lceil\log_2N\rceil
       +L\lceil\log_2(b(w+1))\rceil,\text{ which particularly implies } K2^{-p}\leq N^{-1}.
\]
At any point with $p$ fractional binary digits per coordinate,
layer $j$ can be evaluated exactly with common denominator $2^pq^j$.
The layer values and affine partial sums have magnitude at most
$[b(w+1)]^j$, so the integer numerators have
$\mathcal O(\log_2([b(w+1)]^j) +\log_2(2^pq^j)) = \mathcal O(L\log_2(w+1) +p + L) = \mathcal{O}(p) $ bits \bg{as $N\to \infty$}.
Summing $M\leq N$ outputs and returning their average as an unreduced
fraction adds only $\mathcal O(\log_2N)$ bits.
Using the assumed bounds on $w$ and $L$, all integers in the computation have
$\mathcal O(\log_2N\,\log_2(2+\log_2N))$ bits.

Given a work budget $N$, the algorithm $\cA_{\rm MC}$ draws $M$ independent samples $Z_j\in 2^{-p}\{0,\dots,2^p-1\}$, taking the $p$ binary digits of each coordinate from its random bits, and returns the exact average $M^{-1}\sum_{j=1}^M f(Z_j)$ for the sample count $M$ chosen below.
Each sample requires $dp$ random bits and
$\mathcal O(Lw^2)$ integer operations for exact forward evaluation.
Standard multiple-precision integer arithmetic
\cite[Chapter~1]{BZ:ModernComputerArithmetic2011}
therefore bounds the work per sample, as well as preprocessing and
output, by a polynomial in $1+\log_2N$.
Processing the samples successively requires memory polynomial in
$1+\log_2N$, hence polynomial in the input length, so addresses occupy a bounded number of words in the
model of \Cref{sec:algorithms}.
Choose fixed integers $A,r\geq1$ such that
$H_N:=A(1+\lceil\log_2N\rceil)^r$ bounds the per-sample work
and the combined preprocessing and output cost.
For sufficiently large $N$, set $M:=\lfloor N/(2H_N)\rfloor$.
Then $M\geq N/(4H_N)$ and the total runtime is at most
$(M+1)H_N\leq N$ for every outcome of the random bits.

For the error analysis, let $X_1,\dots,X_M$ be independent and uniformly distributed in $[0,1)^d$.
The first $p$ binary digits of their coordinates are independent unbiased bits, so we may and will assume that the algorithm uses exactly these digits, i.e., 
$Z_j:=2^{-p}\lfloor2^pX_j\rfloor$.
The Lipschitz bound provides $|f(Z_j)-f(X_j)|\leq K2^{-p}\leq N^{-1}$ so that the classical Monte Carlo variance identity for $M^{-1}\sum_{j=1}^Mf(X_j)$ and \bg{the Minkowski inequality reveal}
\[
    \left(\E\left|\cA_{\rm MC}(\Phi,N)-I(f)\right|^2\right)^{1/2}
    \leq \frac{\|f\|_{L^2([0,1]^d)}}{\sqrt M}+\frac1N
    \leq 2\sqrt{\frac{H_N}{N}}\,
          \|f\|_{L^2([0,1]^d)}+\frac1N.
\]
Markov's inequality for the squared error proves
\eqref{eq:random_bit_mc_error} with probability at least $2/3$.
For $0<\gamma<1/2$, absorbing the logarithmic factor into $N^{1/2-\gamma}$ and \bg{using} $N^{-1}\leq N^{-2\gamma}$ proves~\eqref{eq:random_bit_mc_mixed} for sufficiently large~$C_\gamma$.
\end{proof}

The algorithm $\cA_{\rm MC}$ from \Cref{thm:random_bit_mc} attains every exponent $\gamma<1/2$ in~\eqref{eq:random_bit_mc_mixed} on the network classes of \Cref{thm:growing_dim_nn_quadrature_is_hard,thm:fixed_dim_nn_quadrature_is_hard}, whose randomized versions exclude each exponent $\gamma>1/2$ on one of these classes under randomized SETH.
Hence the exponent $1/2$ is optimal.
% The following corollary expresses optimality of the algebraic exponent $1/2$ in the limit $\eps\downarrow0$.
% \begin{corollary}[Optimality of the algebraic exponent]
% \label{cor:random_bit_mc_optimal}
% For either family of network classes in
% \Cref{thm:growing_dim_nn_quadrature_is_hard,thm:fixed_dim_nn_quadrature_is_hard},
% let $\gamma_\eps^\ast$ denote the supremum of the exponents $\gamma>0$ for which a uniform randomized algorithm with runtime at most $N$ satisfies an error bound of the form~\eqref{eq:random_bit_mc_mixed} with some $C_\gamma>0$ and probability at least $2/3$, uniformly over the corresponding class of neural networks for all sufficiently large $N$.
% \fz{In the fixed-dimensional case, take $d=s$ with $s$ fixed; in the growing-dimensional case, $d$ varies over the input dimensions allowed by the theorem.}
% Under randomized SETH,
% \begin{align*}
%     \frac12\leq\gamma_\eps^\ast\leq\frac12+\eps,
%     \qquad
%     \lim_{\eps\downarrow0}\gamma_\eps^\ast=\frac12.
% \end{align*}
% The random-bit Monte Carlo algorithm $\mathcal{A}_{\mathrm{MC}}$ attains the lower endpoint up to logarithmic factors, as in~\eqref{eq:random_bit_mc_error}.
% \end{corollary}
% \begin{proof}
% With $\theta=\theta_\eps$ and sufficiently large $\kappa$,
% \Cref{thm:random_bit_mc} gives~\eqref{eq:random_bit_mc_mixed}
% for every $0<\gamma<1/2$, so $\gamma_\eps^\ast\geq1/2$.
% The corresponding randomized hardness theorem excludes
% $\gamma>1/2+\eps$, so $\gamma_\eps^\ast\leq1/2+\eps$.
% \end{proof}

\section{Differential equations and matrix-vector multiplication}
\label{sec:computational_tasks}

\subsection{Elliptic PDE with neural network right-hand sides}
\label{sec:pde}

We first consider the Poisson equation and use the positivity
of its Green function to turn the integral gap from
Section~\ref{sec:quad} into an $L^2$-gap between the corresponding
solutions. Let $\Omega:=(0,1)^s$ and, for $f\in L^2(\Omega)$, let $u_f\in H^1_0(\Omega)$ denote the unique weak solution of
\begin{align}\label{eq:poisson}
    -\Delta u_f=f
    \quad\text{in }\Omega,
    \qquad
    u_f=0
    \quad\text{on }\partial\Omega.
\end{align}
\begin{lemma}\label{lem:localized_sat_network}
Let $s,k\in\N$ with $k\geq3$ and let $\delta\in\Q_{>0}$.
There exist constants $C_{k,s}>0$, $D_s>0$, and a finite set
$\theta_{s,\delta}^{\rm loc}\subset\Q$ such that the following holds.

For every $n\geq s$ and every
$\alpha\in\widetilde{\mathcal S}_k(n)$, one can construct in
polynomial time a ReLU network
$\Phi_{\alpha}^{\rm loc}$ with input dimension $s$ satisfying $0\leq\mathcal R_{\Phi_{\alpha}^{\rm loc}}\leq1$ and $\mathrm{supp} \mathcal R_{\Phi_{\alpha}^{\rm loc}} \subseteq \left[1/8,7/8\right]^s$. Its width and the number of non-zero entries in every weight matrix are bounded by $C_{k,s}n$, its depth is bounded by $n+D_s$, and all weights and biases belong to
$\theta_{s,\delta}^{\rm loc}$.

If $\alpha$ is unsatisfiable, then $\mathcal R_{\Phi_{\alpha}^{\rm loc}}\equiv0$. If $\alpha$ is satisfiable, then there exists a measurable set
$Q_\alpha^{\rm loc}\subset[1/4,3/4]^s$ such that $|Q_\alpha^{\rm loc}| \geq 2^{-s}(2+\delta)^{-n}$ and $\mathcal R_{\Phi_{\alpha}^{\rm loc}}= 1$ on $Q_\alpha^{\rm loc}$.
\end{lemma}
\begin{proof}
Set $\operatorname{clip}(t) := \sigma(t) - \sigma(t-1)$ and define $P:\R^s\to[0,1]^s$ componentwise by 
\begin{align*}
    P(x)_j
    :=
    \operatorname{clip}\left(2x_j-\frac12\right).
\end{align*}
For $x\in[1/4,3/4]^s$, we have $P(x)=2x - \boldsymbol{1/2}$, so $P$ maps $[1/4,3/4]^s$ bijectively onto $[0,1]^s$. Next define $L(t) := \operatorname{clip}(8t-1)$ and $R(t) := \operatorname{clip}(7-8t)$, as well as
\begin{align*}
    \chi(x)
    :=
    \min_{j=1,\ldots,s}
    \min\{L(x_j),R(x_j)\}.
\end{align*}
Then $0\leq\chi\leq1$, $\chi(x)=1$ for $x\in[1/4,3/4]^s$, and $\chi(x)=0$ for $x\notin[1/8,7/8]^s$. The identity $\min\{a,b\}=a-\sigma(a-b)$ shows that $\chi$ is the realization of a ReLU network whose size depends only on $s$. Let $\Phi_{\alpha,s,\delta}$ be the fixed-dimensional network from
Proposition~\ref{prop:sat_to_nn}, and define
\begin{align*}
    f_\alpha^{\rm loc}(x)
    :=
    \min\left\{
        \mathcal R_{\Phi_{\alpha,s,\delta}}(P(x)),
        \chi(x)
    \right\}.
\end{align*}
Again using the ReLU realization of the minimum, this function is the
realization of a ReLU network $\Phi_\alpha^{\rm loc}$. The asserted network bounds follow from
Proposition~\ref{prop:sat_to_nn}, since all additional networks have
size depending only on $s$. If $\alpha$ is unsatisfiable,
$\mathcal R_{\Phi_{\alpha,s,\delta}}\equiv0$, and hence
$f_\alpha^{\rm loc}\equiv0$. 

Suppose that $\alpha$ is satisfiable. Proposition~\ref{prop:sat_to_nn}
provides a set $Q_\alpha\subset[0,1]^s$ satisfying $|Q_\alpha|\geq(2+\delta)^{-n}$ and $\mathcal R_{\Phi_{\alpha,s,\delta}}=1$ on $Q_\alpha$. Define $Q_\alpha^{\rm loc} := \boldsymbol{1/4} + (1/2)Q_\alpha$. Then $Q_\alpha^{\rm loc}\subset[1/4,3/4]^s$,
$P(Q_\alpha^{\rm loc})=Q_\alpha$, and
$\chi=1$ on $Q_\alpha^{\rm loc}$. Thus, $f_\alpha^{\rm loc}=1$ on $Q_\alpha^{\rm loc}$, and $|Q_\alpha^{\rm loc}| = 2^{-s}|Q_\alpha| \geq 2^{-s}(2+\delta)^{-n}$. Finally, the cutoff implies
$\operatorname{supp}f_\alpha^{\rm loc}
\subseteq[1/8,7/8]^s$.
\end{proof}

\begin{lemma}
\label{lem:poisson_green_lower_bound}
Let $K:=[1/8,7/8]^s$. There exists a constant $c_s^{\rm G}>0$
such that, for every $f\in L^\infty(\Omega)$ with $f\geq0$ and
$\operatorname{supp}f\subseteq K$, the solution of~\eqref{eq:poisson}
satisfies
\begin{align}\label{eq:green_solution_lower_bound}
    \norm{u_f}{L^2(\Omega)}
    \geq
    c_s^{\rm G}\norm{f}{L^1(\Omega)}.
\end{align}
\end{lemma}
\begin{proof}
We use the fundamental solution and the maximum principle to obtain
a lower bound for the Green function when its two arguments are
close. We then partition $K$ into small cubes to obtain the asserted
$L^2$-estimate.

For $s=1$, the Dirichlet Green function on $(0,1)$ is
\[
    G_\Omega(x,y)=\min\{x,y\}\bigl(1-\max\{x,y\}\bigr).
\]
If $x,y\in K=[1/8,7/8]$, both factors are at least $1/8$, so
$G_\Omega(x,y)\geq1/64$. The representation
$u_f(x)=\int_K G_\Omega(x,y)f(y)\,dy$ therefore gives
\[
    \norm{u_f}{L^2(\Omega)}
    \geq\norm{u_f}{L^2(K)}
    \geq\frac{\sqrt{3/4}}{64}\norm{f}{L^1(\Omega)}.
\]
Thus we may take $c_1^{\rm G}:=\sqrt{3/4}/64$.

Now let $s\geq2$, put $r:=\operatorname{dist}(K,\partial\Omega)=1/8$,
and denote the surface area of the unit sphere in $\R^s$ by
$\omega_s$. The fundamental solution of $-\Delta$ is
$\Gamma_s(x):=\varphi_s(|x|_2)$, where
\[
    \varphi_s(\rho):=
    \begin{cases}
        -\dfrac{1}{2\pi}\log\rho,&s=2,\\[1ex]
        \dfrac{\rho^{2-s}}{(s-2)\omega_s},&s\geq3.
    \end{cases}
\]
For fixed $y\in K$, let $W_y\in H^1(\Omega)$ be the harmonic
function with boundary trace $\Gamma_s(\cdot-y)$, that is,
\[
    -\Delta W_y=0\quad\text{in }\Omega,
    \qquad W_y=\Gamma_s(\cdot-y)\quad\text{on }\partial\Omega.
\]
The pole $y$ is separated from the boundary, so these boundary data
are smooth in a neighbourhood of $\partial\Omega$. The Dirichlet
Green function is obtained by subtracting this harmonic correction:
\[
    G_\Omega(x,y)=\Gamma_s(x-y)-W_y(x),\qquad x\neq y.
\]
We use the standard Green-function construction and weak maximum
principle on the cube; see~\cite{giltrud}. Since $\varphi_s$ is decreasing and $|z-y|_2\geq r$ for every
$z\in\partial\Omega$, the maximum principle gives
\[
    W_y(x)\leq\sup_{z\in\partial\Omega}\Gamma_s(z-y)
    \leq\varphi_s(r),\qquad x\in\Omega.
\]
This is an upper bound on $W_y$, not on its absolute value, so it
also applies to the logarithmic fundamental solution in dimension
two without any restriction on $\operatorname{diam}(\Omega)$.
Consequently, whenever $0<|x-y|_2\leq r/2$,
\begin{align*}
    G_\Omega(x,y)
    &\geq\varphi_s(|x-y|_2)-\varphi_s(r)
    \geq\varphi_s(r/2)-\varphi_s(r)
    =:V_s>0,
\end{align*}
with 
\begin{align*}
    V_s&:=
    \begin{cases}
        \dfrac{\log2}{2\pi},&s=2,\\[1ex]
        \dfrac{(2^{s-2}-1)r^{2-s}}{(s-2)\omega_s},&s\geq3.
    \end{cases}
\end{align*}
The Green function is also non-negative for all $x\in\Omega\setminus\{y\}$.
Indeed, for $0<\varepsilon<r/2$ it is harmonic on
$\Omega\setminus\overline{B_\varepsilon(y)}$, has zero trace on
$\partial\Omega$, and is at least $V_s$ on
$\partial B_\varepsilon(y)$. The minimum principle on this punctured
domain gives $G_\Omega(\cdot,y)\geq0$ there. Letting
$\varepsilon\downarrow0$ proves the claim.

For bounded $f$, the Green representation of the weak solution is
\[
    u_f(x)=\int_K G_\Omega(x,y)f(y)\,dy
    \quad\text{for almost every }x\in\Omega.
\]
The singularity at $x=y$ is locally integrable, since
$\int_0^r\rho^{s-1}|\varphi_s(\rho)|\,d\rho<\infty$.
Thus this representation applies to the present $L^\infty$ data;
see~\cite{giltrud}.

The diameter of $K$ is too large to apply the bound $V_s$ to all
pairs in $K\times K$ directly. Set
$q:=\lceil12\sqrt s\rceil$ and partition $K$, up to boundaries of
measure zero, into $M:=q^s$ congruent cubes $Q_1,\ldots,Q_M$.
Each has side length $3/(4q)$, and hence
\[
    \operatorname{diam}(Q_i)=\frac{3\sqrt s}{4q}\leq\frac1{16}
    =\frac r2,
    \qquad |Q_i|=\frac{|K|}{M}.
\]
Put $m_i:=\int_{Q_i}f(y)\,dy$. Non-negativity of $G_\Omega$ and $f$
allows us to discard the contribution from $K\setminus Q_i$.
For almost every $x\in Q_i$, the local lower bound therefore yields
\[
    u_f(x)\geq\int_{Q_i}G_\Omega(x,y)f(y)\,dy\geq V_s m_i.
\]
Integrating the squared bounds over the disjoint cubes and using
$(\sum_{i=1}^M m_i)^2\leq M\sum_{i=1}^M m_i^2$, we obtain
\begin{align*}
    \norm{u_f}{L^2(\Omega)}^2
    &\geq\sum_{i=1}^M\int_{Q_i}u_f(x)^2\,dx
    \geq V_s^2\frac{|K|}{M}\sum_{i=1}^M m_i^2 \geq V_s^2\frac{|K|}{M^2}
    \left(\sum_{i=1}^M m_i\right)^2
    =V_s^2\frac{|K|}{M^2}\norm{f}{L^1(\Omega)}^2.
\end{align*}
Taking square roots proves~\eqref{eq:green_solution_lower_bound}
with the explicit choice
\[
    c_s^{\rm G}:=\frac{V_s(3/4)^{s/2}}{\lceil12\sqrt s\rceil^s}>0,
    \qquad s\geq2.
\]
All quantities in this constant depend only on $s$.
\end{proof}

\begin{corollary}\label{cor:pde_sat_gap}
Let $f_\alpha^{\rm loc}
=\mathcal R_{\Phi_\alpha^{\rm loc}}$ be the function from
Lemma~\ref{lem:localized_sat_network}, let $u_\alpha$ solve
\begin{align*}
    -\Delta u_\alpha=f_\alpha^{\rm loc},
    \qquad
    u_\alpha|_{\partial\Omega}=0,
\end{align*}
and set $U_\alpha:=\norm{u_\alpha}{L^2(\Omega)}$. Then
\begin{align*}
    \alpha\text{ unsatisfiable} \Longrightarrow
    U_\alpha=0, \quad \text{while} \quad
    \alpha\text{ satisfiable} \Longrightarrow
    U_\alpha
    \geq
    c_s^{\rm G}2^{-s}(2+\delta)^{-n}.
\end{align*}
Moreover,
\begin{align}\label{eq:pde_source_solution_scale}
    \norm{f_\alpha^{\rm loc}}{L^2(\Omega)}^2
    \leq
    (c_s^{\rm G})^{-1}U_\alpha.
\end{align}
\end{corollary}
\begin{proof}
If $\alpha$ is unsatisfiable, then
$f_\alpha^{\rm loc}\equiv0$ and hence $u_\alpha=0$. If $\alpha$ is
satisfiable, Lemmas~\ref{lem:localized_sat_network}
and~\ref{lem:poisson_green_lower_bound} give
\begin{align*}
    U_\alpha
    \geq
    c_s^{\rm G}\norm{f_\alpha^{\rm loc}}{L^1(\Omega)}
    \geq
    c_s^{\rm G}|Q_\alpha^{\rm loc}|
    \geq
    c_s^{\rm G}2^{-s}(2+\delta)^{-n}.
\end{align*}
Finally, $0\leq f_\alpha^{\rm loc}\leq1$ and
Lemma~\ref{lem:poisson_green_lower_bound} imply $\norm{f_\alpha^{\rm loc}}{L^2(\Omega)}^2 \leq \norm{f_\alpha^{\rm loc}}{L^1(\Omega)} \leq (c_s^{\rm G})^{-1}U_\alpha$.
\end{proof}

\begin{theorem}\label{thm:poisson_solver_hard}
For every $s\in\N$ there exists $C_s>0$ such that for every
$\eps>0$ there exist $C_\eps>0$ and a finite set
$\theta_\eps\subset\Q$ such that the following holds.

Let $C>0$. Suppose $\cA$ is a deterministic algorithm which, for
every $N\geq2$, accepts every ReLU network $\Phi$ satisfying
\begin{itemize}
    \item input dimension $s$;
    \item width at most $C_\eps(1+\log_2N)$;
    \item depth at most $C_s(1+\log_2N)$;
    \item at most $C_\eps(1+\log_2N)$ non-zero entries in every
    weight matrix;
    \item all weights and biases belong to $\theta_\eps$;
\end{itemize}
and returns $\cA(\Phi,N)\in L^2(\Omega)$ in runtime at most $N$.
Suppose that the output representation allows
$\norm{\cA(\Phi,N)}{L^2(\Omega)}$ to be computed in runtime at most
$N$. Let $u_\Phi\in H^1_0(\Omega)$ solve $-\Delta u_\Phi= \mathcal R_\Phi$. Suppose that, for some $\gamma>0$ and every admissible $(\Phi,N)$,
\begin{align*}
    \norm{
        \cA(\Phi,N)-u_\Phi
    }{L^2(\Omega)}
    \leq
    C\norm{\mathcal R_\Phi}{L^2(\Omega)}N^{-\gamma}
    +
    C^2N^{-2\gamma}.
\end{align*}
Then SETH implies $\gamma\leq 1/2+\eps$.
\end{theorem}
\begin{proof}
Use the localized networks $\Phi_\alpha^{\rm loc}$ from
Lemma~\ref{lem:localized_sat_network} and the work parameter
$N=N_n$ from the proof of
Theorem~\ref{thm:fixed_dim_nn_quadrature_is_hard}.
Set $U_\alpha:=\norm{u_\alpha}{L^2(\Omega)}$. The error bound and a reverse triangle inequality show
\begin{align*}
\left|\norm{\cA(\Phi_\alpha^{\rm loc},N)}{L^2(\Omega)} - U_\alpha \right| \leq C\fz{\norm{\mathcal R_{\Phi_{\alpha}^{\rm loc}}}{L^2(\Omega)}}N^{-\gamma} + C^2N^{-2\gamma}.
\end{align*}
Moreover, \Cref{cor:pde_sat_gap} provides $c_s^{\mathrm{G}}>0$ \fz{(which we may replace by a smaller positive rational constant)} such that
\begin{align*}
    \alpha\text{ unsatisfiable} \Longrightarrow U_\alpha=0, \qquad  \alpha\text{ satisfiable} \Longrightarrow
    U_\alpha\geq q_n:= c_s^{\rm G}2^{-s}(2+\delta)^{-n},
\end{align*}
and~\eqref{eq:pde_source_solution_scale} controls $\norm{f_\alpha^{\rm loc}}{L^2(\Omega)}^2 \leq (c_s^{\rm G})^{-1}U_\alpha$.
The remaining estimates and threshold test in the proof of Theorem~\ref{thm:growing_dim_nn_quadrature_is_hard} apply with $q_n$ as above, $I_\alpha$ replaced by $U_\alpha$, and $C$ replaced by $C\max\{1,(c_s^{\rm G})^{-1/2}\}$ for all sufficiently large $n$.
This concludes the proof.
\end{proof}

\subsection{Matrix-vector multiplication for neural network encoded matrices}
\label{sec:mvm}

A neural network with $2d$ inputs provides a succinct representation
of a $2^d\times2^d$ matrix. More precisely, let
$\Phi$ be a ReLU network with input dimension $2d$ and define
\begin{align*}
    (M_\Phi)_{ij}
    :=
    \mathcal R_\Phi
    \big(
        b(i),b(j)
    \big),
    \qquad
    i,j\in\{0,\ldots,2^d-1\},
\end{align*}
where $b(i)\in\{0,1\}^d$ denotes the binary representation of $i$. Although $M_\Phi$ has exponentially many entries in $d$, the network
$\Phi$ may contain only $\mathcal O(d)$ non-zero parameters. We ask
whether this succinct representation permits matrix-vector
multiplication substantially faster than explicitly forming the
matrix.

\begin{lemma}\label{lem:boolean_matrix_mvm_gap}
Let $D\in\N$ and consider any non-zero Boolean matrix $M\in\{0,1\}^{D\times D}$. Then $\norm{M\boldsymbol 1}{2} \geq\norm{M}{2}\norm{\boldsymbol 1}{2} D^{-1/2}$, where $\boldsymbol1\in\R^D$ is the all-ones vector.
In particular, $\norm{M\boldsymbol1}{2}\geq1$.
\end{lemma}
\begin{proof}
For $i=1,\ldots,D$, let $r_i:=\sum_{j=1}^D M_{ij}$ be the number of non-zero entries in the $i$-th row. Since every $r_i$ is a non-negative integer, $r_i\leq r_i^2$ for every $i$. Furthermore, $M_{ij} \in \{0, 1\}$ and $M_{ij}^2 = M_{ij}$ holds for all $i,j$. Hence,
\begin{align*}
    \norm{M}{2}^2
    \leq
    \norm{M}{F}^2
    =
    \sum_{i=1}^D \sum_{j=1}^D M_{ij}^2 = \sum_{i=1}^D r_i \leq \sum_{i=1}^D r_i^2 = \sum_{i=1}^D \left(\sum_{j=1}^D M_{ij}\right)^2 = \norm{M\boldsymbol1}{2}^2.
\end{align*}
Therefore, $\norm{M}{2} \leq \norm{M\boldsymbol1}{2}$. Since $\norm{\boldsymbol1}{2}=\sqrt D$, we obtain the first assertion. Finally, since $M$ is non-zero, at least one row satisfies $r_i\geq1$,
and hence
$\norm{M\boldsymbol1}{2}\geq1$.
\end{proof}

\begin{lemma}\label{lem:sat_to_boolean_matrix}
Let $k\geq3$, let $d\in\N$, and let
$\alpha\in\widetilde{\mathcal S}_k(2d)$.
Fix, for example, $\delta_0:=1$ and let
$\Psi_{\alpha,\delta_0}$ be the network from
Proposition~\ref{prop:sat_to_nn}.
Then
\begin{align*}
    (M_{\Psi_{\alpha,\delta_0}})_{ij}
    =
    \mathcal I_\alpha\big(b(i),b(j)\big)
    \in\{0,1\}.
\end{align*}
Consequently, $\alpha$ is unsatisfiable if and only if $M_{\Psi_{\alpha,\delta_0}}=0$. If $\alpha$ is satisfiable, then $\norm{M_{\Psi_{\alpha,\delta_0}}\boldsymbol1}{2} \geq1$ and $\norm{M_{\Psi_{\alpha,\delta_0}}}{2} \norm{\boldsymbol1}{2} \leq 2^{d/2} \norm{M_{\Psi_{\alpha,\delta_0}}\boldsymbol1}{2}$.
\end{lemma}

\begin{proof}
The first identity follows, because the neural network agrees exactly
with the Boolean evaluation on the corners of the hypercube. Hence,
$M_{\Psi_{\alpha,\delta_0}}$ is the truth table of $\alpha$, with the
first $d$ variables indexing its rows and the remaining $d$ variables
indexing its columns.

Thus the matrix is zero exactly when $\alpha$ has no satisfying
assignment. If $\alpha$ is satisfiable, the matrix is a non-zero
Boolean matrix. Lemma~\ref{lem:boolean_matrix_mvm_gap}, applied with
$D=2^d$, yields the remaining assertions.
\end{proof}

\begin{theorem}\label{thm:nn_matrix_vector_hard}
For every $\eps>0$ there exists $C_\eps>0$ and a finite set
$\theta_{\rm MVM}\subset\Q$ such that the following holds.

Let $C>0$, and suppose that $\cA$ is a deterministic algorithm which,
for every $N\geq2$, accepts every ReLU network $\Phi$ with an even
input dimension $2d$ satisfying \fz{$2^d\leq N$ and}
\begin{itemize}
    \item $2d\leq C_\eps(1+\log_2N)$;
    \item the width is at most
    $C_\eps(1+\log_2N)$;
    \item the depth is at most $4$;
    \item every weight matrix contains at most
    $C_\eps(1+\log_2N)$ non-zero entries;
    \item every weight and bias belongs to $\theta_{\rm MVM}$,
\end{itemize}
and returns $\cA(\Phi,N)\in\R^{2^d}$ in runtime at most $N$. 
\fz{The condition $2^d\leq N$ allows for the explicit output vector; without a restriction of this kind, the output size alone could make the runtime requirement impossible.}

Suppose that, for some $\gamma>0$ and every admissible $(\Phi,N)$,
\begin{align*}
    \norm{
        \cA(\Phi,N)-M_\Phi\boldsymbol1
    }{2}
    \leq
    C\norm{M_\Phi}{2}
    \norm{\boldsymbol1}{2}N^{-\gamma}
    +
    C^2N^{-2\gamma}.
\end{align*}
Then SETH implies $\gamma\leq1/4+\eps$.
\end{theorem}
\begin{proof}
Fix $\eps>0$ and suppose, for contradiction, that $\gamma > 1/4 + \eps$. Choose a rational number $\omega$ such that
\begin{align*}
    \max\left\{
        \frac12,\frac{1}{1+4\eps}
    \right\}
    <\omega<1.
\end{align*}
By Lemma~\ref{lem:sparse_seth}, SETH implies
$\widetilde s_k\to1$. Hence we can choose $k\geq3$ such that $\omega<\widetilde s_k$. Choose a rational number $B$ satisfying $2^\omega<B<2^{\widetilde s_k}$. We set $\theta_{\rm MVM} := \theta_{{\rm CNF},1}$.

It suffices to consider formulas with an even number of variables.
Indeed, an unused dummy variable may be added to any formula, changing
the number of variables by at most one and preserving satisfiability
and sparsity.

Let $\alpha\in\widetilde{\mathcal S}_k(2d)$ and put $n:=2d$ as well as $N_n:=\left\lceil B^n\right\rceil$. Since $B>2^\omega$ and $\omega>1/2$, we have $2^d=2^{n/2}<B^n\leq N_n$. In particular, the dimension of the output vector is at most the available runtime.

Let $\Psi_\alpha:=\Psi_{\alpha,1}$ be the $n$-dimensional network from
Proposition~\ref{prop:sat_to_nn}.
It has depth $4$, width $\mathcal O_k(n)$, and
$\mathcal O_k(n)$ non-zero entries in each weight matrix.
Moreover, $n \leq {\left(\log_2N_n\right)}/{\left(\log_2B\right)}$. Thus, after choosing $C_\eps$ sufficiently large,
$\Psi_\alpha$ is an admissible input for
$\cA(\,\cdot\,,N_n)$. Write $M_\alpha:=M_{\Psi_\alpha}$ and $y_\alpha:=M_\alpha\boldsymbol1$. By Lemma~\ref{lem:sat_to_boolean_matrix}, if $\alpha$ is unsatisfiable, then $M_\alpha=0$ and $y_\alpha=0$, whereas, if $\alpha$ is satisfiable, then $\norm{y_\alpha}{2}\geq1$ and $\norm{M_\alpha}{2} \norm{\boldsymbol1}{2} \leq 2^{d/2}\norm{y_\alpha}{2} = 2^{n/4}\norm{y_\alpha}{2}$. 

By the assumption $\gamma>1/4+\eps$, we have $\gamma\omega> 1/4$. Since $\log_2B>\omega$, we also have $\gamma\log_2B>1/4$. Therefore, $\rho := 2^{1/4}B^{-\gamma} <1$. Consequently,
\begin{align*}
    C2^{n/4}N_n^{-\gamma}
    \leq
    C2^{n/4}B^{-\gamma n}
    =
    C\rho^n.
\end{align*}
Choose $n_0$ sufficiently large that, for all $n\geq n_0$, we have $C\rho^n\leq 1/8$ and $C^2B^{-2\gamma n}\leq 1/8$. 

If $\alpha$ is unsatisfiable, the approximation error bound becomes $\norm{\cA(\Psi_\alpha,N_n)}{2} \leq C^2N_n^{-2\gamma} \leq 1/8$.
In particular, each output coordinate is trivially bounded by $|(\cA(\Psi_\alpha,N_n))_i|\leq 1/8$ for all $i\in\CSrange{2^d}$.

If $\alpha$ is satisfiable, then $\norm{y_\alpha}{2}\geq1$ provides
\begin{align*}
    \norm{
        \cA(\Psi_\alpha,N_n)-y_\alpha
    }{2} \leq
    C2^{n/4}N_n^{-\gamma}
    \norm{y_\alpha}{2}
    +
    C^2N_n^{-2\gamma} \leq
    \frac18\norm{y_\alpha}{2}
    +
    \frac18\leq \frac14\norm{y_\alpha}{2}.
\end{align*}
Since $M_{\alpha}$ is a Boolean matrix, $y_\alpha$ is nonzero with nonnegative integer entries and at least one output coordinate of $\cA(\Psi_\alpha,N_n)$ must exceed $1/2$.
Indeed, otherwise each coordinate with $(y_\alpha)_i\geq1$ would satisfy
\begin{align*}
    \left|\cA(\Psi_\alpha,N_n)_i-(y_\alpha)_i\right|
    \geq (y_\alpha)_i-\frac12
    \geq \frac12(y_\alpha)_i,
\end{align*}
contradicting the preceding error bound.
Thus, satisfiability can be decided by testing whether any output coordinate exceeds $1/2$.
For $n<n_0$, we use brute force.

The network construction takes polynomial time in $n=\mathcal{O}(\log N_n)$.
Comparing any rational output coordinate $p/q$, with $q>0$, to $1/2$ reduces to the integer comparison $2p>q$ with work linear in the number of words encoding that coordinate.
Summing over all coordinates gives work linear in the total output length and hence $\mathcal O(N_n)$, since writing the output is included in the runtime of $\cA$.
Hence the resulting algorithm for sparse $k$-SAT has runtime $\mathcal O(N_n) = \mathcal O(B^n)$. 
Since $\log_2B<\widetilde s_k$, this contradicts the definition of $\widetilde s_k$.
\end{proof}

\begin{corollary}\label{cor:mvm_eps_four}
\fz{Assume SETH and fix $0<c<1/2$. For every $0<\xi<1$, there is
a constant $C_\xi>0$ for which no uniform deterministic algorithm
has worst-case runtime $\mathcal O(D^{2-\xi})$ on all depth-$4$
networks $\Phi$ with $2\log_2D$ inputs, width and layerwise sparsity
at most $C_\xi\log_2D$, parameters in the fixed finite set
$\theta_{\rm MVM}$, and $M_\Phi\in\{0,1\}^{D\times D}$,
while returning an explicit vector $\widetilde y$ satisfying
\begin{align*}
    \norm{\widetilde y-M_\Phi\boldsymbol1}{2}
    \leq
    cD^{-1/2}\norm{M_\Phi}{2}\norm{\boldsymbol1}{2}.
\end{align*}
Here $D$ ranges over powers of two. For non-zero matrices this is
a relative-error requirement with target $\varepsilon=cD^{-1/2}$;
for the zero matrix it requires exact zero output.
Equivalently, no uniform worst-case bound $\mathcal O(\varepsilon^{-4+\zeta})$
with fixed $\zeta>0$ holds over all these classes.}
\end{corollary}
\begin{proof}
\fz{Suppose such an algorithm exists. Choose $k$ with
$\widetilde s_k>1-\xi/2$ and encode a sparse $k$-SAT formula in
$n=2d$ variables using \Cref{lem:sat_to_boolean_matrix}.
Choose $C_\xi$ large enough to include these networks, and put
$D=2^{n/2}$. In the unsatisfiable case the output must be zero.
Otherwise, \Cref{lem:boolean_matrix_mvm_gap} gives
$\norm{\widetilde y-M_\Phi\boldsymbol1}{2}
\leq c\norm{M_\Phi\boldsymbol1}{2}$.
Since $c<1/2$, the coordinate test from the preceding proof
distinguishes these cases. The resulting SAT algorithm has runtime
$\mathcal O(\operatorname{poly}(n)+2^{(1-\xi/2)n})$,
contradicting the choice of $k$. An unused variable treats odd $n$.
Finally, substituting $\varepsilon=cD^{-1/2}$ turns
$\varepsilon^{-4+\zeta}$ into a constant multiple of
$D^{2-\zeta/2}$, which gives the equivalent formulation.}
\end{proof}

\section{Extensions and limitations of the results}\label{sec:extensions}

\subsection{Uniformly bounded weights}

The finite parameter sets in our main results depend on $\eps$. In
particular, the factor $a_\eta=(1-2\eta)^{-1}$ in
Lemma~\ref{lem:r_with_relu} becomes large when $\eta$ approaches
$1/2$. This is not essential. By using additional layers, the same
transition can be realized with uniformly bounded weights.

\begin{lemma}\label{lem:bounded_weights}
For every $m\geq2$, there exists a ReLU network
$\Psi_m\in\NN(1,2,m+1,\{-2,-1,0,1,2\})$ such that
\[
\cR_{\Psi_m}(t)=
\begin{cases}
0,&t\leq 1/2,\\
2^m(t-1/2),&1/2<t<1/2+2^{-m},\\
1,&t\geq1/2+2^{-m}.
\end{cases}
\]
In particular, if $\eta_m:=1/2-2^{-m}$, then
$\cR_{\Psi_m}=r_{\eta_m}$.
\end{lemma}

\begin{proof}
For $m=1$, the function
$\sigma(2t-1)-\sigma(2t-2)$ has the required form. Assume that
$\Psi_m$ has been constructed and let
$T(t):=\sigma(2t)-\sigma(2t-1)$. Then $T$ is represented by a
depth-two ReLU network with parameters in
$\{-2,-1,0,1,2\}$ and $\cR_{\Psi_{m+1}}:=T\circ\cR_{\Psi_m}$ has the required form with $m$ replaced by $m+1$. Composition
increases the depth by one and preserves width two and the parameter
set.
\end{proof}
\begin{corollary}\label{cor:bounded_weight_hardness}
The conclusions of
Theorems~\ref{thm:growing_dim_nn_quadrature_is_hard}
and~\ref{thm:fixed_dim_nn_quadrature_is_hard} remain valid if the
admissible networks are required to have uniformly bounded weights
and biases.

More precisely, the hard networks can be chosen such that all weights
and biases have absolute value at most $2$. In the growing-dimensional
case the depth is bounded by a constant depending on $\eps$, while in
the fixed-dimensional case it remains bounded by
$C_s(1+\log_2N)$ after increasing the constant.
\end{corollary}

\begin{proof}
In the proofs of
Theorems~\ref{thm:growing_dim_nn_quadrature_is_hard}
and~\ref{thm:fixed_dim_nn_quadrature_is_hard}, choose $m$ sufficiently
large that
$\omega\log_2(2+\delta_m)<\widetilde s_k$, where
$\delta_m=2/(2^{m-1}-1)$, and use $\Psi_m$ from
Lemma~\ref{lem:bounded_weights} in place of the depth-two realization
of $r_\eta$.

The CNF construction uses only copies of this ramp and affine maps
whose coefficients have absolute value at most $2$. Hence all
parameters can be kept uniformly bounded. When composing the building
blocks, their non-negative outputs may be passed through an additional
ReLU layer instead of merging adjacent affine maps; this prevents the
creation of larger product weights. Since $m$ depends only on
$\eps$ and not on the number $n$ of Boolean variables, this changes
the depth only by a constant depending on $\eps$. The width and the
number of non-zero parameters in every layer remain of order
$\mathcal O_k(n)$. The satisfiable and unsatisfiable cases are
unchanged because the realization of the ramp is exactly
$r_{\eta_m}$.
\end{proof}
The factor $\eta_m^{-1}=2+\delta_m$ in the fixed-dimensional encoder
can be written as $2(1+\delta_m/2)$ and therefore implemented by two
successive affine layers whose coefficients have absolute value at
most $2$. The same bounded-weight modification applies to the fixed-domain and
Poisson reductions, since the additional affine maps use only fixed
constants and can be decomposed into a fixed number of affine maps
with coefficients bounded by two.

\subsection{The role of depth}

Our main quadrature lower bounds use networks with more than one hidden
layer. It is natural to ask how much depth is actually needed for the
underlying hardness mechanism. We show that already two hidden layers
are sufficient to encode SAT in the zero set of the integral. The
quantitative integral gap is weaker than in
Proposition~\ref{prop:sat_to_nn}, which explains why this construction
does not give the sharp convergence exponents from Section~\ref{sec:quad}. We also refer to~\cite{np-hard} for a related result on the hardness of computing function norms.

\begin{lemma}\label{lem:shallow_sat}
Let $k\in\N$ and $0<\eta<1/2$. For every
$\alpha\in\mathcal S_k(n)$, one can construct in time polynomial in
the size of $\alpha$ a ReLU network $\Phi_{\alpha,\eta}^{\rm sh}$
with input dimension $n$ and depth $3$ such that
$0\leq\cR_{\Phi_{\alpha,\eta}^{\rm sh}}\leq1$\fz{ on $[0,1]^n$}. Its hidden-layer widths and the number of non-zero entries in every
weight matrix are bounded by $C_{k,\eta}\#\alpha$. Consequently, its
width is at most $\max\{n,C_{k,\eta}\#\alpha\}$, and all weights and
biases belong to a finite set $\theta_{k,\eta}\subset\Q$ if
$\eta\in\Q$. If $\alpha$ is unsatisfiable, then
$\cR_{\Phi_{\alpha,\eta}^{\rm sh}}\equiv0$\fz{ on $[0,1]^n$}. If $\alpha$ is
satisfiable, then there exists a set $Q_\alpha\subseteq[0,1]^n$ with
$|Q_\alpha|\geq(\eta/k)^n$ such that
$\cR_{\Phi_{\alpha,\eta}^{\rm sh}}=1$ on $Q_\alpha$. In particular,
\[
\alpha\text{ satisfiable}
\quad\Longrightarrow\quad
\int_{[0,1]^n}\cR_{\Phi_{\alpha,\eta}^{\rm sh}}(x)\,dx
\geq(\eta/k)^n.
\]
\end{lemma}
\begin{proof}
Fix a clause $C\in\alpha$ and let
$i_1,\ldots,i_q$, $q\leq k$, denote the distinct variables occurring
in $C$. Let $Z_C\subseteq\{0,1\}^q$ be the set of assignments of these
variables which satisfy $C$. For $z\in Z_C$, define
\[
d_{C,z}(x):=
\sum_{j:z_j=0}x_{i_j}
+
\sum_{j:z_j=1}(1-x_{i_j})
\]
and
\[
h_{C,z}(x):=
\frac{\sigma(1/2-d_{C,z}(x))
-\sigma(\eta-d_{C,z}(x))}
{1/2-\eta}.
\]
Since $d_{C,z}$ is affine, $h_{C,z}$ is represented by two neurons in
one hidden layer. We have $0\leq h_{C,z}\leq1$,with 
$h_{C,z}(x)=1$ if $d_{C,z}(x)\leq\eta$, and $h_{C,z}(x)=0$ if $d_{C,z}(x)\geq1/2$.

If $z,z'\in\{0,1\}^q$ are distinct, then
$|z-z'|_1\geq1$. Hence it is impossible to have simultaneously
$d_{C,z}(x)<1/2$ and $d_{C,z'}(x)<1/2$, because otherwise
\begin{align*}
1 &\leq |z-z'|_1 = \sum_{j:z_j\neq z'_j}|z_j-z'_j| = \sum_{j:z_j = 0, z'_j=1} \left(x_{i_j} + (1-x_{i_j})\right) + \sum_{j:z_j = 1, z'_j=0} \left((1-x_{i_j}) + x_{i_j}\right) \\
&\leq \sum_{j:z_j = 0} x_{i_j} + \sum_{j:z_j = 1} (1-x_{i_j}) + \sum_{j:z'_j = 0} x_{i_j} + \sum_{j:z'_j = 1} (1-x_{i_j}) = d_{C,z}(x) + d_{C,z'}(x) < 1,
\end{align*}
Thus, for
$G_C(x):=\sum_{z\in Z_C}h_{C,z}(x)$ we have
$0\leq G_C(x)\leq1$. We now define
\[
F_\alpha(x):=
\sigma\left(1+\sum_{C\in\alpha}(G_C(x)-1)\right).
\]
This is the realization of a network with two hidden layers. To avoid
a bias depending on $\#\alpha$, we include one constant neuron
$\sigma(1)=1$ in the first hidden layer for every clause and connect
it to the second hidden layer with weight $-1$. Hence all parameters
belong to a finite set depending only on $k$ and $\eta$. Since every
clause contains at most $k$ variables and $|Z_C|\leq2^k$, the
hidden-layer widths and layerwise sparsity are bounded by
$C_{k,\eta}\#\alpha$.

Let now $x\in[0,1]^n$, and let $z(x)\in\{0,1\}^n$ be obtained by
rounding every coordinate of $x$ to a closest Boolean value. If
$\alpha$ is unsatisfiable, there exists a clause $C_x$ which is false
under $z(x)$. Every assignment $z\in Z_{C_x}$ therefore differs from
$z(x)$ in at least one variable occurring in $C_x$. For this variable
we have $|x_i-z_i|\geq1/2$, and hence
$d_{C_x,z}(x)\geq1/2$. Thus $G_{C_x}(x)=0$. Since
$G_C(x)\leq1$ for every other clause,
\[
1+\sum_{C\in\alpha}(G_C(x)-1)\leq0,
\]
and therefore $F_\alpha(x)=0$. Since $x$ was arbitrary,
$F_\alpha\equiv0$.

Suppose now that $\alpha$ is satisfiable and let
$z_\alpha\in\{0,1\}^n$ be a satisfying assignment. Consider
\[
Q_\alpha:=
\left\{x\in[0,1]^n:
|x_i-z_{\alpha,i}|\leq\eta/k
\text{ for all }i\right\}.
\]
For every clause $C$ the restriction of $z_\alpha$ to the variables
of $C$ belongs to $Z_C$, and for every $x\in Q_\alpha$ its distance
from $x$ is at most $k(\eta/k)=\eta$. Hence $G_C(x)=1$ for every
clause and therefore $F_\alpha(x)=1$. Finally,
$|Q_\alpha|=(\eta/k)^n$.
\end{proof}

The shallow construction does not recover the sharp quantitative
lower bound from Theorem~\ref{thm:growing_dim_nn_quadrature_is_hard}.
Indeed, its satisfiable integral is only bounded below by
$(\eta/k)^n$. For fixed $k$ the same SETH argument would give the
conditional bound
\[
\gamma\leq
\frac{\log_2(k/\eta)}{2\widetilde s_k}
\]
for the \fz{mixed error bound} considered in
Theorem~\ref{thm:growing_dim_nn_quadrature_is_hard}, provided
$\widetilde s_k>0$. Since SETH obtains
$\widetilde s_k\to1$ only by letting $k\to\infty$, this bound becomes
worse rather than better as $k$ grows. We therefore use the shallow
construction only to establish the qualitative complexity barrier
above.

\subsection{Quasi-Monte Carlo quadrature}
The recent result~\cite{dick2026} shows that randomly initialized neural networks have favorable properties for QMC quadrature, i.e., suitable QMC points improve upon the convergence rates bounds in this work. 
In the following, we demonstrate empirically, that QMC does not work robustly for the integrands used in the hardness constructions above.

Let $C\subset\{0,1\}^d$ and define
\begin{align}\label{eq:rquad}
    f_C(x)
    :=
    \sum_{z\in C}
    \min_{1\leq j\leq d}
    r_\eta\big(1-|z_j-x_j|\big).
\end{align}
For the experiments we choose $\eta=0.2$ and a uniformly sampled
subset $C\subset\{0,1\}^d$ with $|C|=2^{d-1}$.

For every $z\in C$, the corresponding term is equal to one in a
neighborhood of $z$ and vanishes once at least one coordinate is
sufficiently far from $z$. Thus, $f_C$ has a structure similar to the
localized Boolean regions appearing in our SAT construction.

\begin{figure}[tbp]
    \centering
    \begin{minipage}{0.49\textwidth}
    \includegraphics[width=\textwidth]{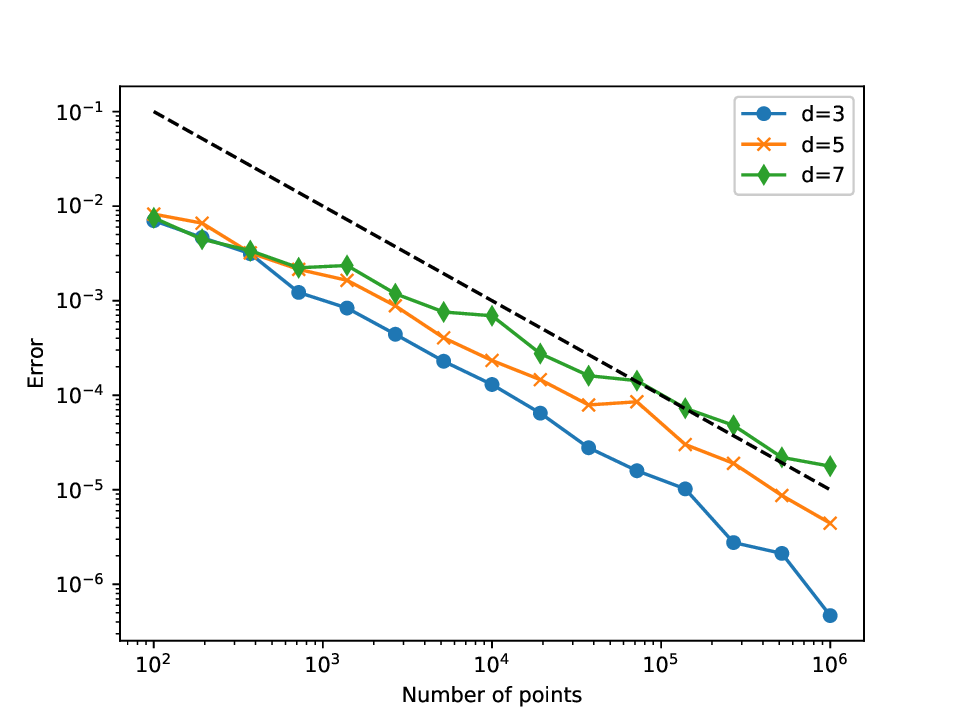}
    \end{minipage}
    \begin{minipage}{0.49\textwidth}
    \includegraphics[trim=52 52 52 52,clip,width=\textwidth]{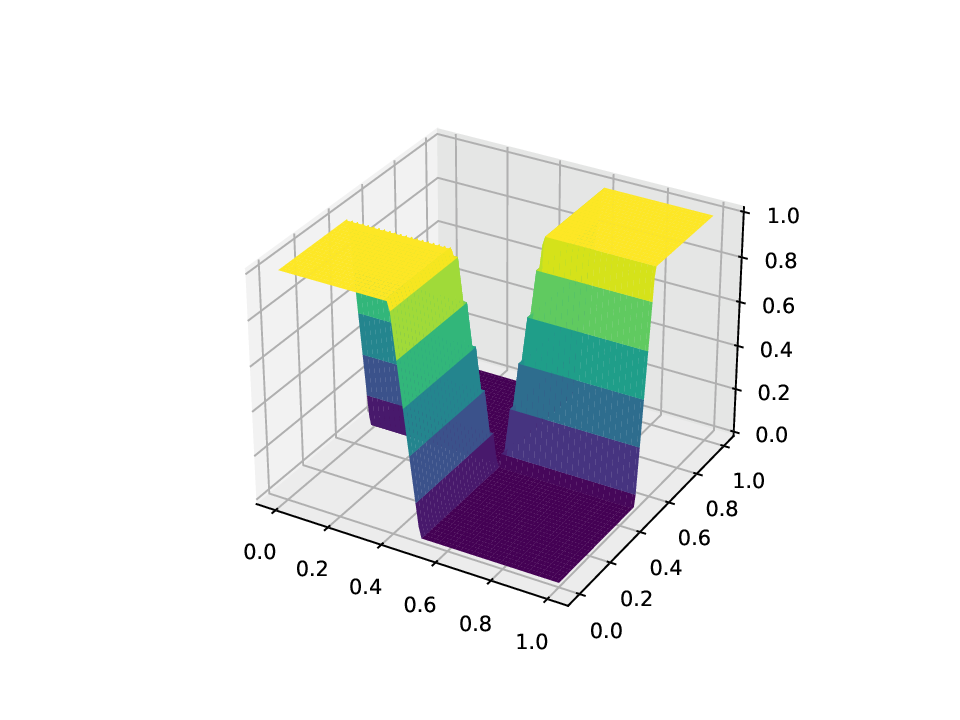}
    \end{minipage}
    \caption{Left: Absolute Sobol quadrature error for the
    hard-instance-inspired integrands from~\eqref{eq:rquad}; the
    dashed line has slope $-1$. Right: An example of such an integrand
    in dimension two.}
    \label{fig:rquad}
\end{figure}

The results are shown in Figure~\ref{fig:rquad} suggest that the worst-case instances
used in the complexity reduction are substantially less benign than
typical randomly initialized or trained networks. Identifying
structural assumptions on neural networks under which QMC admits
provably faster convergence remains an interesting open problem.

\FloatBarrier
\subsection{Fast quadrature for one-hidden-layer networks}

The hardness results above require more than one hidden layer. In
fixed dimension, ridge functions with a ReLU nonlinearity can be
integrated explicitly, and so can networks with one hidden
layer.

\begin{proposition}\label{prop:relu_ball_integral}
Let $R>0$ and let \(B_R(c)\subset\R^d\) be equipped with normalized Lebesgue
measure. For \(a \in\R^d\) and \(b\in\R\), set $A:=R|a|_2$ and $B := b- a\cdot c$. If \(A=0\), then
\[
    \int_{B_R(c)}\sigma(a \cdot x-b)\,d\mu(x) = \sigma(-B).
\]
If \(A>0\), let \(u:=B/A\). Then
\begin{align*}
\int_{B_R(c)}\sigma(a \cdot x-b)\,d\mu(x) =
\begin{cases}
-B,
    & B \leq -A,\\[1mm]
\displaystyle
\frac{A(1-u^2)^{(d+1)/2}}{(d+1)\beta(1/2,(d+1)/2)} - \frac{B}{2}
\left(
1-\operatorname{sgn}(u)
I_{u^2}\left(\frac12,\frac{d+1}{2}\right)
\right),
    & |B|<A,\\[3mm]
0,
    & B \geq A,
\end{cases}
\end{align*}
where \(\beta\) is the beta function and \(I_x\) is the regularized
incomplete beta function.
\end{proposition}
\begin{proof}
After the change of variables \(x=c+Ry\) and a rotation of
coordinates, the integral equals
\[
    \frac{1}{\beta(1/2,(d+1)/2)}
    \int_{-1}^1
    \sigma(At-B)(1-t^2)^{(d-1)/2}\,dt.
\]
The cases \(B\leq-A\) and \(B\geq A\) are immediate. If
\(|B|<A\), integration starts at \(u=B/A\). The identities
\[
    \int_u^1 t(1-t^2)^{(d-1)/2}\,dt
    =
    \frac{(1-u^2)^{(d+1)/2}}{d+1}
\]
and
\[
    \frac{1}{\beta(1/2,(d+1)/2)}
    \int_u^1(1-t^2)^{(d-1)/2}\,dt
    =
    \frac12\left(
        1-\operatorname{sgn}(u)
        I_{u^2}\left(\frac12,\frac{d+1}{2}\right)
    \right)
\]
give the stated formula.
\end{proof}

On the hypercube, the algorithm is a bit more complicated.
\begin{proposition}\label{prop:relu_cube_fourier}
Let $d\geq2$, $\ba\in(\R\setminus\{0\})^d$, and $b\in\R$. Set $\be:=\ba/|\ba|_2$, $t:=b/|\ba|_2$, and $L:=\sum_{j=1}^d|e_j|$.
If $t\geq L$, then
$\int_{[-1,1]^d}\sigma(\ba\cdot x-b)\,dx=0$, while for $t\leq-L$
the integral equals $-2^d b$. If $|t|<L$ we define
\[
    P_{\be}(u):=
    \prod_{j=1}^d\operatorname{sinc}(e_j u),
    \qquad
    \operatorname{sinc}(v):=
    \begin{cases}
        \sin(v)/v,&v\neq0,\\
        1,&v=0.
    \end{cases}
\]
Then
\begin{align}\label{eq:relu_cube_fourier}
\int_{[-1,1]^d}\sigma(\ba\cdot x-b)\,dx =
2^{d-1}\left(|\ba|_2 2\sqrt{d} - b\right) +
\frac{2^d |\ba|_2}{\pi}
\int_0^\infty
P_{\be}(u)
\frac{\cos\left(u2\sqrt{d}\right)-\cos(tu)}{u^2}\,du.
\end{align}
At $u=0$, the integrand is defined by its limit $(t^2- 4d)/2$.
Each evaluation of the one-dimensional integrand costs
$\mathcal O(d)$ arithmetic operations and elementary function evaluations.
\end{proposition}
\begin{proof}
The two cases $t\geq L$ and $t\leq-L$ follow from
$-|\ba|_2L\leq\ba\cdot x\leq |\ba|_2L$ on the cube and
$\int_{[-1,1]^d}\ba\cdot x\,dx=0$. Hence, suppose $|t|<L$.
Complete $\be$ to an orthonormal basis
$\be,\bv_2,\ldots,\bv_d$ of $\R^d$. Writing
$x=s\be+s_2\bv_2+\cdots+s_d\bv_d$ gives $\ba\cdot x=|\ba|_2s$.
For each $s$, the remaining coordinates describe the section
\[
    H_{\be,s}\cap[-1,1]^d,
    \qquad H_{\be,s}:=\{x\in\R^d:\be\cdot x=s\}.
\]
Consequently, integration in these coordinates yields
\[
\int_{[-1,1]^d}\sigma(\ba\cdot x-b)\,dx
=\int_t^{2\sqrt{d}}(|\ba|_2s-b)|H_{\be,s}\cap[-1,1]^d|_{d-1}\,ds.
\]
Here $|\cdot|_{d-1}$ denotes the $(d-1)$-dimensional volume.
The upper limit $2\sqrt d$ is sufficient because the section is
empty for $|s|>L$, and $L\leq\sqrt d$ by the Cauchy--Schwarz inequality. The section volumes admit the Fourier representation
\cite{intersections,intersections2}
\begin{align*}
|H_{\be,s}\cap[-1,1]^d|_{d-1}
&=
\frac{2^{d-1}}{\pi}
\int_{-\infty}^{\infty}
P_{\be}(u)\cos(su)\,du.
\end{align*}
Since every $e_j$ is non-zero, $P_{\be}(u)=\mathcal O(|u|^{-d})$
as $|u|\to\infty$. Thus $d\geq2$ ensures absolute integrability,
and we may interchange the integrals in $s$ and $u$. This gives
\begin{align*}
\int_{[-1,1]^d}\sigma(\ba\cdot x-b)\,dx
&=\frac{2^{d-1}}{\pi}\int_{-\infty}^{\infty}
P_{\be}(u)\int_t^{2\sqrt{d}}(|\ba|_2s-b)\cos(su)\,ds\,du.
\end{align*}
Using $|\ba|_2t=b$, direct integration of the inner integral gives
\[
\int_t^{2\sqrt{d}}(|\ba|_2s-b)\cos(su)\,ds
=
(|\ba|_2 2\sqrt{d} -b)\frac{\sin(u2\sqrt{d})}{u}
+|\ba|_2\frac{\cos(u2\sqrt{d})-\cos(tu)}{u^2}.
\]
In particular, the integrand in the outer integral can be written as
\begin{align*}
f_{\ba,b}(u)
&:=P_{\be}(u)\left(
(|\ba|_2 2\sqrt{d} -b)\frac{\sin(u2\sqrt{d})}{u}
-2|\ba|_2\frac{\sin((2\sqrt{d}+t)u/2)\sin((2\sqrt{d}-t)u/2)}{u^2}
\right).
\end{align*}
This is an even function, so the integral over $\R$ is twice the
integral over $[0,\infty)$. Both terms are constant multiples of products of sinc functions.
Such integrals are often called Borwein integrals. In particular,
for $\boldsymbol\beta\in\R^r$, the identity
\[
\int_0^\infty\prod_{j=1}^r\operatorname{sinc}(\beta_j u)\,du
=\frac{\pi}{2|\beta_k|}
\qquad\text{if }|\beta_k|>\sum_{j\neq k}|\beta_j|
\]
is available from~\cite{borwein}. For the first term of $f_{\ba,b}$,
we may apply it with $\boldsymbol\beta=(e_1,\ldots,e_d,2\sqrt{d})$, because
$2\sqrt d>L$. Therefore,
\[
\int_0^\infty P_{\be}(u)\frac{\sin(u2\sqrt{d})}{u}\,du=\frac\pi2.
\]
Evaluating this term explicitly concludes the proof.
\end{proof}

For numerical quadrature, it is useful to write the remaining integrand as
\begin{align*}
h_{\be,t}(u) := P_{\be}(u)\frac{\cos(u2\sqrt{d})-\cos(tu)}{u^2} =-\frac{4d -t^2}{2}\,P_{\be}(u)
\operatorname{sinc}\left(\frac{2\sqrt{d}+t}{2}u\right)
\operatorname{sinc}\left(\frac{2\sqrt{d}-t}{2}u\right).
\end{align*}
This expression avoids subtracting two nearly equal cosines near
$u=0$ and shows directly that $h_{\be,t}$ is analytic there. To approximate the Borwein integrals, choose $T,p\in\N$,
truncate the integral at $T$, and use $p$-point Gauss--Legendre
quadrature on each interval $[i,i+1]$, $i=0,\ldots,T-1$.
Products of sinc functions are analytic and bounded on every strip
of fixed width around the real axis. Standard Gaussian quadrature
estimates therefore give an error $\mathcal O(Te^{-\kappa p})$
on $[0,T]$ for some $\kappa>0$; see~\cite{quadrature}.
Moreover, $h_{\be,t}(u)=\mathcal O(u^{-d-2})$ as $u\to\infty$,
so the omitted tail has size $\mathcal O(T^{-d-1})$.
Together, the two errors yield
\[
    \mathcal O\left(Te^{-\kappa p}+T^{-d-1}\right)
    \qquad\text{at a cost of }\mathcal O(dTp)
\]
arithmetic operations and elementary function evaluations.
The constants depend on the fixed parameters $\ba$, $b$, and $d$;
in particular, the tail bound deteriorates when a non-zero $e_j$
approaches zero. Here $d$ is the dimension after removing zero
coordinates. Finally, integrals over $[0,1]^d$ follow by changing variables.

\printbibliography

\end{document}